\documentclass{amsart}
\usepackage{amsmath,amsthm,bbm}
\usepackage{amsfonts,amssymb}

\newtheorem{thm}{Theorem}[section]
\newtheorem{cor}[thm]{Corollary}
\newtheorem{lem}[thm]{Lemma}
\newtheorem{prop}[thm]{Proposition}

\newtheorem{defn}[thm]{Definition}
\newtheorem{re}[thm]{Remark}

\theoremstyle{remark}

\def\CH{{\mathcal H}}

\def\CN{{\mathcal N}}

\def\CY{{\mathcal Y}}
\def\CC{{\mathbb C}}

\def\NN{{\mathbb N}}
\def\PP{{\mathbb P}}

\def\RR{{\mathbb R}}

\def\supp{\operatorname{supp}}

\newcommand{\wh}{\widehat}

\def\cS{{\mathcal S}}

\def\cX{{\mathcal X}}

\newcommand{\eps}{{\varepsilon}}

\def\R{{\mathbb R}}

\def\a{{\alpha}}
\def\b{{\beta}}
\def\g{{\gamma}}
\def\ha{{\wh a}}
\def\ha{{\widehat a}}
\def\hb{{\widehat b}}

\def\Ph{{\mathcal H}}
\def\PP{{\mathcal H}}
\def\ph{{\varphi}}
\def\cH{{\mathcal H}}
\def\cS{{\mathcal S}}

\def\HH{{\bf H}}
\def\cM{{\mathcal M}}
\def\ONE{{\mathbbm 1}}

\def\L2{{L^2}}
\def\Linfty{{\infty}}

\def\wtd{\widetilde}

\def\Fapq{{\F^{\a p}_q}}
\def\bb{b}
\def\bapq{{\b^{\a p}_q}}
\def\Bapq{{\B^{\a p}_q}}

\def\V{V}
\def\W{W}
\def\LL{\Lambda}

\def\T{D}
\def\HH{\cH}

\def\Fapq{{F^{\a q}_p}}

\def\Bapq{{B^{\a q}_p}}
\def\bapq{{b^{\a q}_p}}

\def\cd{c_\diamond}

\def\u{u}
\def\v{v}

\def\hcX{\widehat{\cX}}

\def\dm{\diamond}
\def\up{\lambda}
\def\mm{\omega}

\def\del{\lambda}

\begin{document}

\title{Decomposition of spaces of distributions induced by Hermite expansions}

\author{Pencho Petrushev and Yuan Xu}
\address{Department of Mathematics\\University of South Carolina\\
Columbia, SC 29208.}
\email{pencho@math.sc.edu}
\address{Department of Mathematics\\ University of Oregon\\
Eugene, Oregon 97403-1222.}\email{yuan@math.uoregon.edu}


\subjclass{ 42B35, 42C15} \keywords{Localized kernels, frames,
Hermite polynomials, Triebel-Lizorkin spaces, Besov spaces}
\thanks{The second author was partially supported by the
NSF under Grant DMS-0604056.}

\begin{abstract}
Decomposition systems with rapidly decaying elements (needlets) based on
Hermite functions are introduced and  explored. It is proved that
the Triebel-Lizorkin and Besov spaces on $\R^d$ induced by Hermite expansions
can be characterized in terms of the needlet coefficients.
It is also shown that the Hermite Triebel-Lizorkin and Besov spaces
are, in general, different from the respective classical spaces.
\end{abstract}

\maketitle

\pagestyle{myheadings}
\thispagestyle{plain}
\markboth{PENCHO PETRUSHEV AND YUAN XU}
         {DECOMPOSITION OF SPACES INDUCED BY HERMITE EXPENSIONS}

\section{Introduction}\label{introduction}
\setcounter{equation}{0}

The purpose of this paper is to extend the fundamental results of
Frazier and Jawerth \cite{F-J1, F-J2}
on the $\varphi$-transform to the case of Hermite expansions on $\R^d$.
In the spirit of \cite{F-J1, F-J2} we will construct a pair of dual frames
in terms of Hermite functions and use them to characterize
the Hermite-Triebel-Lizorkin and Hermite-Besov spaces.

Let $\{h_n\}_{n=0}^\infty$ be the $L^2(\R^d)$ normalized univariate Hermite functions
(see \S\ref{Local-kernels}). The $d$-dimensional Hermite functions are defined by
$
\HH_\alpha(x) := h_{\alpha_1}(x_1)\cdots h_{\alpha_d}(x_d).
$
Then the kernel of the orthogonal projector of $L^2$ onto
$\W_n:= {\rm span}\, \{\HH_\alpha: |\alpha|=n\}$ is given by
$ 
\HH_n(x,y):= \sum_{|\alpha| = n} \HH_\alpha(x)\HH_\alpha(y).
$ 
Our construction of Hermite frames hinges on the fundamental fact that
for compactly supported $C^\infty$ functions $\ha$ the kernels
$ 
\LL_n(x,y):= \sum_{j=0}^\infty \ha(\frac{j}{n}) \HH_j(x,y)
$ 
decay rapidly away from the main diagonal in $\R^d$.
This fact was established in \cite{Epp2} for dimension $d=1$
and in \cite{Dzub} in general.
We obtain a more precise estimate in
Theorem~\ref{thm:derivative} below.
We utilize kernels of such kind for the construction of a pair of dual frames
$\{\ph_\xi\}_{\xi\in\cX}$, $\{\psi_\xi\}_{\xi\in\cX}$,
where $\cX$ is a multilevel index set.
The frame elements have almost exponential localization
(see (\ref{local-Needlets21})) which prompted us to call them ``needlets".
The needlet systems of this article can be viewed as an analogue of the $\varphi$-transform
of Frazier and Jawerth \cite{F-J1, F-J2}.
Frames of the same nature in the case $d=1$ have been previously introduced
in \cite{Epp2}.

Our primary goal is to utilize needlets to the characterization of
the Triebel-Lizorkin and Besov spaces in the context of Hermite expansions.
To be more specific, assume that $\ha\in C^\infty$,
$\supp \ha \subset [1/4, 4]$, and
$|\ha|>c$ on $[1/3, 3]$,
and define
\begin{equation}\label{def-Phi-Intr}
\Phi_0 := \CH_0
\quad\mbox{and}\quad
\Phi_j := \sum_{\nu=0}^\infty \ha
\Big(\frac{\nu}{4^{j-1}}\Big)\PP_\nu, \quad j\ge 1.
\end{equation}
Then for all appropriate indices we define the Hermite-Triebel-Lizorkin space
$F^{\alpha q}_p=F^{\alpha q}_p(H)$ as the set of all tempered distributions
$f$ such that
$$
\|f\|_{F^{\alpha q}_p} :=
\Big\|\Big(\sum_{j=0}^\infty
(2^{\alpha j}|\Phi_j*f(\cdot)|)^q\Big)^{1/q}\Big\|_p
< \infty,
$$
where $\Phi_j*f(x):=\langle f, \overline{\Phi(x, \cdot)} \rangle$
(see Definition~\ref{def.Tr-Liz}).
We define the Hermite-Besov spaces $B^{\alpha q}_p=B^{\alpha q}_p(H)$ by the norm
$$
\|f\|_{B^{\alpha q}_p} :=
\Big(\sum_{j=0}^\infty \Big(2^{\alpha j}\|\Phi_j*f\|_p\Big)^q\Big)^{1/q}.
$$
One normally uses binary dilations in (\ref{def-Phi-Intr})
(see e.g. \cite[\S 10.3]{T1} and also \cite{Dzub, Epp1}).
We dilate $\ha$ by factors of $4^j$ instead since then the Hermite F- and B-spaces
embed just as the classical F- and B-spaces.

Our main results assert that the Hermite-Triebel-Lizorkin and Hermite Besov
spaces can be characterized in terms of respective sequence norms of
the needlet coefficients of the distributions
(Theorems~\ref{thm:needlet-Tr-Liz}, \ref{thm:character-Besov}).
Furthermore, we use these results to show that the Hermite-F- and B-spaces
of essentially positive smoothness are different from the respective
classical F- and B-spaces on $\R^d$.

Our development here is a part of a bigger project for needlet characterization
of Triebel-Lizorkin and Besov spaces on nonclassical domains such as
the unit sphere \cite{NPW}, the interval with Jacobi weights \cite{KPX1},
and the unit ball \cite{KPX2}.


The rest of the paper is organized as follows:
Section~\ref{preliminaries} contains some background material.
The needlets are introduced in \S\ref{def-needlets}.
In \S\ref{Tri-Liz-spaces} the Hermite-Triebel-Lizorlin spaces are defined
and characterized via needlets.
The Hermite-Besov spaces are introduced and characterized in \S\ref{Besov-spaces}.
Section~\ref{proofs} contains the proofs of a number of lemmas and theorems
from \S\ref{preliminaries}-\S\ref{Besov-spaces}.

Some useful notation:
$\|f\|_p:= \|f\|_{L^p(\R^d)}$;
for a measurable set $E\subset \R^d$,
$|E|$ denotes the Lebesgue measure of $E$ and
$\ONE_E$ is the characteristic function of $E$.
Also, for $x\in\RR^d$, $|x|$ is the Euclidean norm of $x$,
$|x|_\infty:=\max_{1\le j\le d} |x_j|$, and
$d(x, E):= \inf_{y\in E} |x-y|_\infty$
is the $\ell^\infty$ distance of $x$ from $E\subset \RR^d$.
Positive constants are denoted by $c$, $c_1, \dots$ and they
may vary at every occurrence; $A\sim B$ means $c_1A\le B\le c_2 A$.


\section{Preliminaries}\label{preliminaries}
\setcounter{equation}{0}

\subsection{Localized kernels induced by Hermite functions}\label{Local-kernels}

We begin with a review of some basic properties of Hermite polynomials
and functions. (For background information we refer the reader to \cite{Th}.)
The Hermite polynomials are defined by
$$
H_n(t) = (-1)^n e^{t^2} \Big(\frac{d}{dt}\Big)^n \Big(e^{-t^2}\Big),
\qquad n = 0, 1, \ldots .
$$
These polynomials are
orthogonal with respect to $e^{-t^2}$ on $\RR$.
We will denote the $L^2$-normalized Hermite functions by
$$
h_n(t) := \left(2^n n! \sqrt{\pi}\right)^{-1/2} H_n(t) e^{-t^2/2}.
$$
One has
$$
\int_{\RR} h_n(t) h_m(t) dt = \left(2^n n! \sqrt{\pi} \right)^{-1}
\int_{\RR} H_n(t) H_m(t) e^{-t^2} dt = \delta_{n,m}.
$$
As is well known the Hermite functions form an orthonormal basis for $\L2(\RR)$.

As already mentioned, the $d$-dimensional Hermite functions $\HH_\alpha$
are defined by
\begin{equation}\label{def-HH}
\HH_\alpha(x) := h_{\alpha_1}(x_1) \cdots h_{\alpha_d}(x_d),
\quad \alpha=(\alpha_1, \dots, \alpha_d).
\end{equation}
Evidently $e^{|x|^2/2} \HH_\alpha(x)$ is a polynomial of degree
$|\alpha|:= \alpha_1+\cdots+ \alpha_d$. The Hermite functions form an
orthonormal basis for $\L2(\RR^d)$.
Moreover, $\HH_\alpha$ are eigenfunctions of
the Hermite operator $\T:= - \Delta + |x|^2$ and
\begin{equation}\label{Herm-oper}
\T \HH_\alpha = (2 |\alpha|+d)\HH_\alpha,
\end{equation}
where $\Delta$ is the Laplacian.
The operator $\T$ can be written in the form
\begin{equation}\label{rep-T}
\T = \frac{1}{2} \sum_{j=1}^d (A_j A_j^* + A_j^* A_j),
\quad\hbox{where}\quad A_j= -\frac{\partial}{\partial x_j}+x_j,
\quad A_j^* = \frac{\partial}{\partial x_j}+x_j.
\end{equation}
Let $e_j $ denote the $j$th coordinate vector in $\RR^d$.
Then the operators $A_j$ and $A_j^*$ satisfy
\begin{equation}\label{eq:Aj}
A_j\HH_\alpha = (2\alpha_j +2)^{\frac{1}{2}} \HH_{\alpha+e_j}
\quad\mbox{and}\quad
A_j^*\HH_\alpha = (2\alpha_j)^{\frac{1}{2}} \HH_{\alpha-e_j}.
\end{equation}
Combining these two relations shows that $\{\HH_\alpha\}$ satisfy
the recurrence relation
\begin{equation}\label{eq:xPhi}
x_j \HH_\alpha(x)
= \left(\tfrac{\alpha_j +1}{2}\right)^{\frac{1}{2}} \HH_{\alpha+e_j}(x)
+ \left(\tfrac{\alpha_j}{2}\right)^{\frac{1}{2}} \HH_{\alpha-e_j}(x)
\end{equation}
and also
\begin{equation}\label{eq:DPhi}
\frac{\partial}{\partial x_j} \HH_\alpha(x)
= -\left(\tfrac{\alpha_j +1}{2}\right)^{\frac{1}{2}} \HH_{\alpha+e_j}(x)
+ \left(\tfrac{\alpha_j}{2}\right)^{\frac{1}{2}} \HH_{\alpha-e_j}(x).
\end{equation}

\smallskip

Let $\W_n :=\rm{span}\,\{\HH_\alpha: |\alpha|=n\}$ and
 $\V_n := \bigoplus_{j=0}^n \W_j$.
The kernels of orthogonal projectors on $\W_n$ and $\V_n$ are given by
\begin{equation}\label{def-Kn}
\HH_n(x,y):= \sum_{|\alpha| = n} \HH_\alpha(x)\HH_\alpha(y)
\quad \hbox{and} \quad
K_n (x,y): = \sum_{j=0}^n \HH_j(x,y),
\end{equation}
respectively.

An important role will be played by operators whose kernels are
obtained by smoothing out the coefficients of the kernel $K_n$ by sampling
a compactly supported $C^\infty$ function $\wh a$.
For our purposes we will be considering ``smoothing" functions $\ha$
that satisfy:


\begin{defn} \label{defn:admissible}
A function $\ha \in C^\infty[0, \infty)$
is said to be admissible of type

$(a)$ if
$\supp \ha \subset [0,1+\v]$ $(v>0)$ and $\ha(t)=1$ on $[0, 1]$,
%
and of type

$(b)$ if
$\supp \ha \subset [\u,1+\v]$, where $0<\u<1$, $v>0$.
\end{defn}

\noindent
For an admissible function $\wh a$ we consider the kernel
\begin{equation} \label{def-LLn}
\LL_n(x,y):= \sum_{j=0}^\infty \ha\Big(\frac{j}{n}\Big) \HH_j(x,y).
\end{equation}

It will be critical for our further development that
the kernels $\LL_n(x,y)$ and their derivatives decay rapidly
away from the main diagonal $y = x$ in $\RR^d \times \RR^d$:


\begin{thm} \label{thm:derivative}
Suppose $\wh a$ is admissible in the sense of Definition~\ref{defn:admissible}
and let $\alpha \in \NN_0^d$.
Then for any $k\ge 1$ there exists
a constant $c_k$ depending only on $k$, $\alpha$, $d$, and $\wh a$
such that
\begin{equation}\label{eq:derivative1}
\Big| \frac{\partial^\alpha}{\partial x^\alpha} \LL_n(x,y) \Big|
\le c_k \frac{n^{\frac{|\alpha|}{2}}
[K_{n+[vn]+|\alpha|+k}(x,x)]^{\frac12}
[K_{n+[vn]+k}(y,y)]^{\frac12} } {(1+n^\frac{1}{2}|x-y|)^k}.
\end{equation}
Here the dependence of $c_k$ on $\ha$ is of the form
$c_k=c(k, |\alpha|, \u, d)\max_{0\le l \le k}\|\ha^{(l)}\|_{\Linfty}$.
\end{thm}

We relegate the somewhat lengthy proof of this theorem to \S\ref{proof-loc-est}.

The function
\begin{equation} \label{def:lambda-n}
\lambda_n(x) := \frac{1}{K_n(x,x)}
\end{equation}
is termed {\em Christoffel function} and it is known
(see e.g. \cite{LL}) to have the following asymptotic in dimension $d=1$:
\begin{equation} \label{asympt-lambda-n}
\lambda_n(x)
\sim  n^{-1/2} \left(\max \left\{ n^{-2/3}, 1-\frac{|x|}
{\sqrt{2n}}\right\}\right)^{-1/2}
\end{equation}
uniformly for $n\ge 1$ and $x \le \sqrt{2n}(1+ c'n^{-2/3})$,
where $c'>0$ is any fixed constant.
Consequently, for $d=1$ we have
\begin{equation} \label{est:Kn1}
K_n(x,x)
\sim  n^{1/2} \left(\max \left\{ n^{-2/3}, 1-\frac{|x|}
{\sqrt{2n}}\right\}\right)^{1/2},
\quad x \le \sqrt{2n}(1+ c'n^{-2/3}).
\end{equation}

For $d\ge 2$ one has (see \cite[p. 70]{Th})
\begin{equation}\label{est-Hn}
|\HH_n(x,x)| \le c n^{d/2 -1}, \quad x \in \RR^d.
\end{equation}
This along with (\ref{def-Kn}) leads to
\begin{equation}\label{est:Kn-d}
K_n(x, x)\le cn^{d/2}, \quad x \in \RR^d, \quad d\ge 1.
\end{equation}
On the other hand, it is well known that (see e.g. \cite[p. 26]{Th})
\begin{equation}\label{est-hn}
|h_n(x)| \le ce^{-\gamma x^2}, \quad |x| \ge (4n+2)^{1/2}, \quad \gamma>0,
\end{equation}
and
$\|h_n\|_\infty \le cn^{-1/12}$,
which readily imply
\begin{equation}\label{est:Kn-d2}
K_n(x, x)\le ce^{-\gamma' |x|_\infty^2},
\quad\mbox{if \: $|x|_{\infty}:=\max_{1\le j\le d}|x_j|\ge (4n+2)^{1/2}$,}
\end{equation}
where $\gamma'>0$ depends only on $d$.


Now, combining (\ref{eq:derivative1}) with (\ref{est:Kn-d}) and (\ref{est:Kn-d2})
(setting $\gamma^*:=\gamma'/2$)
we arrive at


\begin{cor}\label{cor:derivative}
Under the hypothesis of Theorem~\ref{thm:derivative} we have
\begin{equation}\label{est:derivative2}
\Big| \frac{\partial^\alpha}{\partial x^\alpha} \LL_n(x,y) \Big|
\le c_k \frac{n^{\frac{|\alpha|+d}{2}}} {(1+n^\frac{1}{2}|x-y|)^k},
\quad x \in \RR^d,
\end{equation}
\begin{equation}\label{est:derivative3}
\Big| \frac{\partial^\alpha}{\partial x^\alpha} \LL_n(x,y) \Big|
\le c_k \frac{e^{-\gamma^* |x|^2_\infty}} {(1+n^\frac{1}{2}|x-y|)^k},
\; \mbox{if}\; |x|_{\infty}\ge (4(n+[vn]+|\alpha|+k)+2)^{1/2},
\end{equation}
and
\begin{equation}\label{est:derivative4}
\Big| \frac{\partial^\alpha}{\partial x^\alpha} \LL_n(x,y) \Big|
\le c_k \frac{e^{-\gamma^* |y|^2_\infty}} {(1+n^\frac{1}{2}|x-y|)^k},
\; \mbox{if}\;\; |y|_{\infty}\ge (4(n+[vn]+k)+2)^{1/2}.
\end{equation}
\end{cor}

Note that an estimate similar to (\ref{est:derivative2}) is proved in \cite{Epp2}
when $d=1$ and $\alpha=0$ and in the general case in \cite{Dzub}.
Estimate (\ref{eq:derivative1}) is new.

We now turn to a lower bound estimate.


\begin{thm}\label{thm:lower-bound}
Let $\ha$ be admissible in the sense of Definition~\ref{defn:admissible}
and $|\ha(t)|>c_*>0$ on $[1, 1+\tau]$, $\tau>0$.
Then for any $\eps>0$,
$$
 \int_{\RR^d} |\LL_n(x,y)|^2 dy  \ge c \, n^{d/2} \quad
     \hbox{for}\quad    |x| \le (1-\eps)\sqrt{2(1+\tau)n},
$$
where $c>0$ depends only on $\tau$, $\eps$, $c_*$, and $d$.
\end{thm}

This theorem provides a lower bound for the range where estimate
(\ref{est:derivative2}) (with $\alpha=0$) is sharp.
To indicate the dependence of $K_n$ on $d$, we write $K_{n,d} = K_n$.
Theorem~\ref{thm:lower-bound} is an immediate consequence of (\ref{est:Kn1})
and the following lemma.


\begin{lem}\label{lem:lower-bound}
If $0<\del <1$, $0<\rho<1$,  and $d\ge 1$, then there exists a constant
$c > 0$ such that for $n\ge 2/\del$
\begin{equation} \label{lowerbd}
\sum_{m=[(1-\del)n]}^{n}
\CH_m^2(x,x) \ge c n^{\frac{d-1}{2}} K_{[\rho n], 1}(t,t) \quad
\hbox{if $t:=|x| \le 2 \sqrt{2n+1} $ }.
\end{equation}
\end{lem}

The proof of this lemma is given in \S\ref{proof-loc-est}.


\subsection{Norm relation}
\label{comparison}

For future use we give here the well known relation between different
norms of functions from $V_n$ (see e.g. \cite{LL}):
For $0 < p, q \le \infty$
\begin{equation}\label{norm-relation}
\|g\|_p \le cn^{\frac{d}{2}|1/q-1/p|}\|g\|_q
\quad \mbox{for} \quad g \in \V_n,
\end{equation}
with $c>0$ depending only on $p$, $q$, and $d$.

This estimate can be proved by means of the kernels from (\ref{def-LLn})
with $\ha$ admissible of type (a).

\subsection{Cubature formula}

In order to define our frame elements, we need a cubature formula exact for
products $fg$ with $f, g\in \V_n$.
Such a formula, however, is readily
available using the Gaussian quadrature formula.


\begin{prop}\label{prop:gaussian}\cite{Sz}
Denote by $t_{\nu,n}$, $\nu=1, 2, \dots, n$, the zeros of the Hermite polynomial
$H_{n}(t)$.
The Gaussian quadrature formula
\begin{equation}\label{eq:gaussian}
\int_{\RR} f(t) e^{-t^2} dt \sim \sum_{\nu=1}^n w_{\nu,n}
f(t_{\nu,n}),
\quad w_{\nu,n}:=\lambda_n(t_{\nu,n}) e^{-t_{\nu,n}^2},
\end{equation}
is exact for all polynomials of degree $2n-1$.
Here
$\lambda_n(\cdot)$ is the Christoffel function defined in
$(\ref{def:lambda-n})$.
\end{prop}

The product nature of $e^{-|x|^2}$ enables us to obtain the desired
cubature formula on $\RR^d$ right away.


\begin{prop}\label{prop:cubature}
Let
$\xi_{\alpha,n} := (t_{\a_1,n}, \ldots, t_{\a_d,n})$ and
$\lambda_{\alpha,n} := \prod_{\nu=1}^d \lambda_n(t_{\a_\nu,n})$.
The~cubature formula
\begin{equation}\label{eq:cubature}
\int_{\RR^d} f(x)g(x) dx \sim \sum_{\alpha_1=1}^n \cdots
\sum_{\alpha_d=1}^n \lambda_{\alpha,n} f(\xi_{\alpha,n}) g(\xi_{\alpha,n})
\end{equation}
is exact for all $f\in \V_{\ell}$, $g \in \V_{m}$ with $\ell+m\le 2n-1$.
\end{prop}


We next record some well known properties of the zeros of Hermite polynomials.
Suppose  $\{\xi_\nu\}$ are the zeros of $H_n(t)$ (with $n$ even) ordered so that
\begin{equation}\label{est-zeros}
\xi_{-\frac{n}{2}}<\dots < \xi_{-1}<0<\xi_1 < \cdots < \xi_{\frac{n}{2}},
\quad \xi_{-\nu}=\xi_\nu.
\end{equation}
From \cite{LL} we have
$\xi_{\frac{n}{2}} \le \sqrt{2n+1} - n^{-1/6}$
and uniformly for $|\nu| \le n/2-1$
\begin{equation}\label{est-diff}
\xi_{\nu+1} - \xi_{\nu-1} \sim n^{-1/2}
\left(\max \left\{ n^{-2/3}, 1-\frac{|\xi_{\nu}|}{\sqrt{2n}} \right\} \right)^{-1/2}.
\end{equation}
Consequently, on account of (\ref{asympt-lambda-n})
\begin{equation}\label{weight1}
\lambda_n(\xi_{\nu})\sim \xi_{\nu-1} - \xi_{\nu+1},
\quad |\nu| \le n/2-1 \quad (\xi_0:=0).
\end{equation}

By \cite[(6.31.19)]{Sz}
\begin{equation}\label{zeros}
\frac{\pi(\nu-\frac{1}{2})}{(2n+1)^{1/2}} < \xi_\nu
< \frac{4\nu+3}{(2n+1)^{1/2}},
\quad \nu=1, \dots, n/2.
\end{equation}
From this and (\ref{est-diff}) we have, for any $\eps>0$,
\begin{equation}\label{est-diff2}
\xi_{\nu+1} - \xi_{\nu-1} \sim n^{-1/2}
\quad \mbox{if \: $|\nu| \le (1/2-\eps)n$,}
\end{equation}
and
\begin{equation}\label{est-diff3}
c_1n^{-1/2}\le \xi_{\nu} - \xi_{\nu-1} \le c_2n^{-1/6}
\quad \mbox{if \: $(1/2-\eps)n< |\nu| \le n/2$.}
\end{equation}
Here the constants depend on $\eps$.

It also follows by (\ref{est-diff}) that
\begin{equation}\label{almost-eq}
\xi_{\nu+1} - \xi_{\nu-1} \sim \xi_{\nu} - \xi_{\nu-2},
\quad -n/2+2 \le \nu \le n/2-1.
\end{equation}


For the construction of our frames in Section~\ref{def-needlets} we need
the cubature formulae from Proposition~\ref{prop:cubature}
with
\begin{equation}\label{def-Nj}
n=2N_j, \quad
\mbox{where} \quad
N_j:=[(1+11\delta)(4/\pi)^24^j]+3
\end{equation}
and $0<\delta<1/37$
is an arbitrary (but fixed) constant.

Given $j\ge 0$, let as above
$\xi_\nu$, $\nu= \pm 1, \dots, \pm N_j$,
be the zeros of $H_{2N_j}(t)$.
Let $\cX_j$ be the set of all nodes of cubature $(\ref{eq:cubature})$ with
$n=2N_j$, i.e. $\cX_j$ is the set of all points $\xi_\a:=(\xi_{\a_1}, \dots, \xi_{\a_d})$,
where $0<|\a_\nu| \le N_j$.
Also, for $\xi=\xi_{\a}$ we denote briefly $\lambda_\xi:=\lambda_{\a, N_j}$.
Note that
$\# \cX_j=(2N_j)^d \sim 4^{jd}$.

An immediate consequence of Proposition~\ref{prop:cubature} is the following


\begin{cor}\label{cor:cubature}
The cubature formula
\begin{equation}\label{cubature2}
\int_{\RR^d} f(x)g(x) dx
\sim \sum_{\xi\in \cX_j} \lambda_\xi f(\xi) g(\xi),
\quad
\lambda_{\xi} := \prod_{\nu=1}^d \lambda_{2N_j}(\xi_{\a_\nu}),
\end{equation}
is exact for all
$f\in \V_{\ell}$, $g \in \V_{m}$ with $\ell+m\le 4N_j-1$.
\end{cor}

For later use we now introduce {\it tiles} $\{R_\xi\}$ induced by the points of $\cX_j$.
Set
\begin{align*}
&I_{1}:=[0, (\xi_{1}+\xi_{2})/2], \; I_{-1}:=-I_1,
\\
&I_{\nu}:=[(\xi_{\nu-1}+\xi_{\nu})/2, (\xi_{\nu}+\xi_{\nu+1})/2],
\quad \nu=\pm 2, \dots, \pm N_{j-1},
\quad\mbox{and}\\
& I_{N_j}:=[(\xi_{N_{j-1}}+\xi_{N_j})/2, \xi_{N_j}+2^{-j/6}], \;
I_{-N_j}:= -I_{N_j}.
\end{align*}
For each $\xi=\xi_\a=(\xi_{\a_1}, \dots, \xi_{\a_d})$ in $\cX_j$ we set
\begin{equation}\label{def.Q-xi}
R_\xi:= I_{\a_1}\times I_{\a_2}\times \cdots \times I_{\a_d},
\end{equation}
and also
\begin{equation}\label{def.Q-j}
Q_j:=[\xi_{-N_j}-2^{-j/6}, \xi_{N_j}+2^{-j/6}]^d
=\cup_{\xi\in\cX_j} R_\xi.
\end{equation}
Thus we have associated to each $\xi\in\cX_j$ ($j\ge 0$) a tile $R_\xi$
so that different tiles do not overlap (have disjoint interiors)
and they cover the cube $Q_j\sim[-2^j, 2^j]^d$.

Observe that by the construction of the tiles $\{R_\xi\}$
and (\ref{weight1}) we have
\begin{equation}\label{tiles}
\lambda_\xi \sim |R_\xi|, \quad \xi\in\cX_j.
\end{equation}


By (\ref{est-diff2})
$|R_\xi| \sim 2^{-jd}$ if $\xi=\xi_\alpha\in \cX_j$ with
$|\alpha|_\infty\le (1/2-\delta/2)2N_j= (1-\delta)N_j$.
Assume that $|\xi_\alpha| \le (1+4\delta)2^{j+1}$.
By (\ref{zeros})
$$
|\xi_\alpha|_\infty>\frac{\pi(|\alpha|_\infty-1/2)}{(4N_j+1)^{1/2}}
\quad \mbox{and hence}\quad
\frac{\pi(|\alpha|_\infty-1/2)}{(4N_j+1)^{1/2}}
< (1+4\delta)2^{j+1}.
$$
Using the definition of $N_j$ in (\ref{def-Nj})
it is easy to show that the above inequality implies
$|\alpha|_\infty \le (1-\delta)N_j$.
Consequently, for $\xi\in \cX_j$,
\begin{equation}\label{tiles1}
R_\xi \sim \xi + [-2^{-j}, 2^{-j}]^d
\quad\mbox{and }\quad
|R_\xi| \sim 2^{-jd} \quad
\mbox{if \; $|\xi|_\infty \le (1+4\delta)2^{j+1}$.}
\end{equation}
On the other hand, by (\ref{est-diff2})-(\ref{est-diff3}) it follows that,
in general,
\begin{equation}\label{tiles2}
\xi + [-c_12^{-j}, c_12^{-j}]^d \subset R_\xi
\subset \xi + [-c_22^{-j/3}, c_22^{-j/3}]^d,
\quad \xi\in \cX_j,
\end{equation}
and hence
\begin{equation}\label{tiles3}
c'2^{-jd} \le |R_\xi| \le c''2^{-jd/3}.
\end{equation}

Finally, note that since the zeros of $H_n$ and $H_{n+1}$ interlace,
each $R_\eta\in\cX_{j+\ell}$, $\ell \ge 1$, may intersect at most
finitely many (depending only on $d$) tiles $R_\xi$, $\xi\in\cX_j$.

\subsection{Maximal operator}\label{max-operator}

Let $\cM_s$ be the maximal operator, defined by
\begin{equation}\label{maxOperator}
\cM_s f(x) := \sup_{Q:\, x \in Q}
\left(\frac{1}{|Q|} \int_Q |f(y)|^s\,dy \right)^{1/s},
\quad x\in \RR^d,
\end{equation}
where the $\sup$ is over all cubes $Q$ in $\RR^d$ with sides parallel
to the coordinate axes which contain $x$.

We will need the Fefferman-Stein vector-valued maximal inequality
(see \cite{Stein}):
If $0<p<\infty$, $0<q\le \infty$, and $0 < s < \min \{p, q\}$,
then for any sequence of functions  $f_1, f_2, \dots$ on $\RR^d$
\begin{equation}\label{max-ineq}
\Big\|\Big( \sum_{j=1}^{\infty}
     \left[\cM_s f_j(\cdot)\right]^q \Big)^{1/q}\Big\|_p
\le c \Big\|\Big( \sum_{j=1}^{\infty}|f_j(\cdot)|^q \Big)^{1/q}\Big\|_p,
\end{equation}
where $c = c(p,q,s,d)$.


\subsection{Distributions on \boldmath $\RR^d$}

As is customary, we will denote by $\cS$ the Schwartz class of all functions
$\phi\in C^\infty(\RR^d)$ such that
\begin{equation}\label{schwartz}
P_{\b,\g}(\phi):=\sup_{x}|x^\g D^\b \phi(x)| <\infty
\quad\mbox{for all }\g, \b.
\end{equation}
The topology on $\cS$ is defined by the semi-norms $P_{\b,\g}$.
Then the space $\cS'$ of all temperate distributions is defined
as the set of all continuous linear functionals on $\cS$.
The pairing of $f\in \cS'$ and $\phi\in\cS$ will be denoted by
$\langle f, \phi \rangle := f(\overline{\phi})$ which is consistent
with the inner product
$\langle f, g \rangle := \int_{\RR^d}f \overline{g}dx$ in $\L2(\RR^d)$.

As~a~convenient notation we introduce the following
``convolution":


\begin{defn}\label{def:convolution}
For functions $\Phi: \RR^d\times\RR^d \to \CC$ and $f: \RR^d \to \CC$,
we write
\begin{equation}\label{convolution}
\Phi*f(x) := \int_{\RR^d} \Phi(x, y)f(y)\,dy.
\end{equation}
More generally,
assuming that $f \in \cS'$ and $\Phi: \RR^d\times\RR^d\to\CC$ is such that
$\Phi(x, y)$ belongs to $\cS$ as a function of $y$ $(\Phi(x, \cdot)\in \cS)$,
we define $\Phi*f$ by
\begin{equation}\label{convolution1}
\Phi*f (x) := \langle f, \overline{\Phi(x, \cdot)} \rangle,
\end{equation}
where on the right $f$ acts on $\overline{\Phi(x, y)}$ as a function of $y$.
\end{defn}

We next record some properties of the above ``convolution" that are
well known and easy to prove.


\begin{lem}\label{lem:convolution}
$(a)$ If $f \in \cS'$ and $\Phi (\cdot, \cdot) \in \cS(\RR^d\times \RR^d)$, then
$\Phi*f \in \cS$.
Furthermore $\cH_n*f\in \V_n$.

$(b)$ If $f \in \cS'$, $\Phi (\cdot, \cdot) \in \cS(\RR^d\times \RR^d)$, and
$\phi \in \cS$, then
$\langle \Phi*f, \phi \rangle = \langle f, \overline{\Phi}*\phi \rangle$.

$(c)$ If $f \in \cS'$,
$\Phi (\cdot, \cdot), \Psi (\cdot, \cdot) \in \cS(\RR^d\times \RR^d)$,
and
$\Phi (y, x)= \Phi (x, y)$, $\Psi (y, x)= \Psi (x, y)$,
then
\begin{equation}\label{Psi*Phi*f}
\Psi*\overline{\Phi}*f(x)
= \langle \Psi(x, \cdot), \Phi(\cdot, \cdot)\rangle*f.
\end{equation}
\end{lem}

Evidently the Hermite functions $\{\CH_\a\}$ belong to
the space of test functions~$\cS$.
More importantly the functions in $\cS$ can be
characterized by the coefficients in their Hermite expansions.
Denote
\begin{equation}\label{def-P*}
P_r^{*}(\phi) := \sum_{n=0}^\infty (n+1)^r\|\PP_n*\phi\|_{2}
=\sum_{n=0}^\infty (n+1)^r
\Big(\sum_{|\alpha|=n}|\langle \phi, \CH_\a\rangle|^2\Big)^{1/2},
\quad r\ge 0.
\end{equation}


\begin{lem}\label{lem:char-S} We have
\begin{equation}\label{char-S}
\phi\in\cS \quad \Longleftrightarrow \quad
|\langle \phi, \CH_\a\rangle| \le c_k(|\a|+1)^{-k}
\quad\mbox{for all $\a$ and all $k$.}
\end{equation}
Moreover, the topology in $\cS$ can be equivalently defined by
the semi-norms $P_r^{*}$ from above.
\end{lem}

\noindent
{\bf Proof.}
(a) Assume first that the right-hand side estimates in (\ref{char-S}) hold.
Applying repeatedly identities (\ref{eq:xPhi})-(\ref{eq:DPhi}) one
easily derives the estimate
\begin{equation}\label{char-S1}
\sup_x |x^\g D^\b\cH_\a(x)| \le c(|\a|+1)^{(|\g|+|\b|)/2}
\max_{|\omega|\le |\a|+|\b|+|\g|} \|\cH_\omega\|_{\infty}
\end{equation}
for all indices $\b$ and $\g$,
which implies $\phi\in\cS$ using that
$\phi=\sum_{n=0}^\infty
      \sum_{|\alpha|=n}\langle \phi, \CH_\a\rangle\CH_\alpha$ in $L^2$.

(b) Suppose $\phi\in\cS$.
Using (\ref{Herm-oper}) we have
\begin{align*}
\langle \phi, \cH_\a \rangle
&=
\int_{\RR^d}\cH_\a(x)\phi(x)dx
= \frac{1}{2|\a|+d}\int_{\RR^d}(-\Delta+|x|^2)\cH_\a(x)\phi(x)dx\notag\\
&= \frac{1}{2|\a|+d}\int_{\RR^d}\cH_\a(x) (-\Delta\phi+|x|^2\phi(x))dx,
\end{align*}
where for the last equality we used integration by parts.
Repeating the above procedure $k$ times we obtain a representation
for $\langle \phi, \cH_\a \rangle$ of the form
\begin{equation}\label{char-S2}
\langle \phi, \cH_\a \rangle
= \frac{1}{(2|\a|+d)^k}\int_{\RR^d}\cH_\a(x)
\sum_{|\b|\le 2k, \, |\g|\le 2k}C_{\b,\g} x^\g D^\b\phi(x)dx
\end{equation}
which yields the right-hand side estimates in (\ref{char-S}).

The equivalence of the topologies in $\cS$ induced by the semi-norms
from (\ref{schwartz}) and (\ref{def-P*}) follows easily by
(\ref{char-S1}) and (\ref{char-S2}).
$\qed$


\section{Construction of building blocks (Needlets)}\label{def-needlets}
\setcounter{equation}{0}

We utilize the localized kernels from Theorem~\ref{thm:derivative}
and the cubature formula from Corollary~\ref{cor:cubature} to the construction of
a pair of dual frames consisting of localized functions on~$\RR^d$.

Let $\ha$, $\hb$
satisfy the conditions:  
\begin{equation}\label{ab1}
\ha, \hb \in C^\infty(\R),
\quad\supp \ha, \supp \hb \subset [1/4, 4],
\end{equation}
\begin{equation}\label{ab3}
|\ha(t)|, |\hb(t)|>c>0 \quad \mbox{if}~~
t \in [1/3, 3],
\end{equation}
\begin{equation}\label{ab5}
\overline{\ha(t)}\;\hb(t) + \overline{\ha(4t)}\;\hb(4t) =1
\quad\mbox{if}~~ t \in [1/4, 1].
\end{equation}
Consequently,
\begin{equation}\label{unity1}
\sum_{\nu=0}^\infty \overline{\ha(4^{-\nu}t)}\;\hb(4^{-\nu}t) = 1,
\quad t \in [1, \infty).
\end{equation}

It is easy to see that (see e.g. \cite{F-J2})
if $\ha$ satisfies (\ref{ab1})-(\ref{ab3}), then
there exists $\hb$ satisfying (\ref{ab1})-(\ref{ab3})
such that $(\ref{ab5})$ holds true.

\smallskip

Assuming that $\ha$, $\hb$ satisfy
(\ref{ab1})-(\ref{ab5}),
we define
\begin{align}
&\Phi_0 :=\PP_0,
\quad
\Phi_j := \sum_{\nu=0}^\infty
\ha\Big(\frac{\nu}{4^{j-1}}\Big)\PP_\nu,
\quad j\ge 1, \quad \mbox{and}  \label{def.Phi-j}\\
&\Psi_0 :=\PP_0,
\quad
\Psi_j := \sum_{\nu=0}^\infty
\hb\Big(\frac{\nu}{4^{j-1}}\Big)\PP_\nu,
\quad j\ge 1. \label{def-Psi-j}
\end{align}
%
%
Let $\cX_j$ be the set of the nodes of cubature formula (\ref{cubature2})
from Corollary~\ref{cor:cubature}
and let $\lambda_\xi$ be the coefficients of that cubature formula.
We now define the $j$th level {\em needlets} by
\begin{equation}\label{def-needlets1}
\ph_\xi(x) := \lambda_\xi^{1/2}\Phi_j(x, \xi)
\quad\mbox{and}\quad
\psi_\xi(x) := \lambda_\xi^{1/2}\Psi_j(x, \xi),
\qquad \xi \in \cX_j.
\end{equation}
Write $\cX := \cup_{j = 0}^\infty \cX_j$,
where equal points from different levels $\cX_j$ are considered
as distinct elements of $\cX$.
We use $\cX$ as an index set to define a pair of dual needlet systems
$\Phi$ and $\Psi$ by
\begin{equation}\label{def-needlets2}
\Phi:=\{\ph_\xi\}_{\xi\in\cX}, \quad \Psi:=\{\psi_\xi\}_{\xi\in\cX}.
\end{equation}
According to their further roles, we will call
$\{\ph_\xi\}$ {\em analysis needlets}
and $\{\psi_\xi\}$ {\em synthesis needlets}.

The almost exponential localization of the needlets will be critical
for our further development.
Indeed, by (\ref{est:derivative2}) we have
\begin{equation}\label{local-Needlets2}
|\Phi_j(\xi, x)|, |\Psi_j(\xi, x)|
\le \frac{c_k2^{jd}}{(1+2^{j}|x-\xi|)^k},
\quad x\in \RR^d, \quad \forall k,
\end{equation}
Fix $L>0$. Then by (\ref{est:derivative4}) it follows that
for any $k>0$
\begin{equation}\label{local-Needlets3}
|\Phi_j(\xi, x)|, |\Psi_j(\xi, x)|
\le \frac{c_k2^{-jL}}{(1+2^{j}|x-\xi|)^k},
\quad x\in \RR^d, \;\mbox{if}\;\;
|\xi|_\infty > (1+\delta)2^{j+1}.
\end{equation}
Here $c_k$ depends on $L$ and $\delta$ as well.

From above and (\ref{tiles})-(\ref{tiles2})
we infer
\begin{equation}\label{local-Needlets21}
|\ph_\xi(x)|, |\psi_\xi(x)|
\le \frac{c_k2^{jd/2}}{(1+2^{j}|x-\xi|)^k},
\quad \mbox{if}\;\;
|\xi|_\infty \le (1+\delta)2^{j+1},
\end{equation}
and
\begin{equation}\label{local-Needlets22}
|\ph_\xi(x)|, |\psi_\xi(x)|
\le \frac{c_k2^{-jL}}{(1+2^{j}|x-\xi|)^k},
\quad \mbox{if}\;\;
|\xi|_\infty > (1+\delta)2^{j+1}.
\end{equation}

The following proposition provides a discrete
decomposition of $\cS'$ and $L^p(\RR^d)$ via needlets.


\begin{prop}\label{prop:needlet-rep}
$(a)$ If $f \in \cS'$, then
\begin{equation}\label{Needle-rep}
f = \sum_{j=0}^\infty
\Psi_j*\overline{\Phi}_j*f
\quad\mbox{in} \;\; \cS' \;\; \mbox{and}
\end{equation}
\begin{equation}\label{needlet-rep}
f = \sum_{\xi \in \cX}
\langle f, \ph_\xi\rangle \psi_\xi
\quad\mbox{in} \;\; \cS'.
\end{equation}

$(b)$ If $f \in L^p$, $1\le p \le \infty$, then
$(\ref{Needle-rep})-(\ref{needlet-rep})$ hold in $L^p$.
Moreover, if $1 < p < \infty$, then the convergence in
$(\ref{Needle-rep})-(\ref{needlet-rep})$ is unconditional.
\end{prop}

\noindent
{\bf Proof.}
(a)
By the definition of $\Phi_j$ and $\Psi_j$ in (\ref{def.Phi-j})-(\ref{def-Psi-j})
it follows that
$\Psi_0*\overline{\Phi}_0=\PP_0$
and
$$
\Psi_j*\overline{\Phi}_j(x, y)
=\sum_{\nu=4^{j-2}}^{4^j} \overline{\ha\Big(\frac{\nu}{4^{j-1}}\Big)}
\hb\Big(\frac{\nu}{4^{j-1}}\Big)\PP_\nu(x, y),
\quad j\ge 1.
$$
Note that
$\Psi_j(x, y)$ and $\Phi_j(x, y)$
are symmetric functions (e.g. $\Psi_j(y, x)=\Psi_j(x, y)$)
since $\cH_\nu(x, y)$ are symmetric and hence
$\Psi_j*\overline{\Phi}_j(x, y)$ is well defined.
Now, (\ref{unity1}) and Lemma~\ref{lem:char-S} yield (\ref{Needle-rep}).

To establish (\ref{needlet-rep}), we note that
$\Psi_j(x, \cdot)$ and $\overline{\Phi_j(y, \cdot)}$
belong to $\V_{4^j}$ and
applying the cubature formula from Corollary~\ref{cor:cubature},
we obtain
\begin{eqnarray*}
\Psi_j*\overline{\Phi}_j(x, y)
&=& \int_{\RR^d} \Psi_j(x, u)\overline{\Phi_j(y, u)}\,dy\\
&=&\sum_{\xi\in \cX_j}
\lambda_\xi\Psi_j(x, \xi)\overline{\Phi_j(y, \xi)}
= \sum_{\xi\in \cX_j}\psi_\xi(x)\overline{\ph_\xi(y)}.
\end{eqnarray*}
Consequently,
$$
\Psi_j*\overline{\Phi}_j*f = \sum_{\xi\in \cX_j}
\langle f,\ph_\xi\rangle \psi_\xi.
$$
This along with (\ref{Needle-rep}) implies (\ref{needlet-rep}).

(b) Representation (\ref{Needle-rep}) in $L^p$ follows easily by the rapid decay
of the kernels of the $n$th partial sums. We omit the details.
Then (\ref{needlet-rep}) in $L^p$ follows as above.
The unconditional convergence in $L^p$, $1<p<\infty$, follows by
Proposition~\ref{prop:identification} and Theorem~\ref{thm:needlet-Tr-Liz} below.
\qed


\begin{re}\label{rem:frame}
It is well known that there exists a function $\ha \ge 0$ satisfying
$(\ref{ab1})-(\ref{ab3})$ such that
$\ha^2(t) + \ha^2(4t) =1$, $t \in [1/4, 1]$.
Suppose that in the above construction
$\hb = \ha$ and $\ha\ge 0$.
Then
$\ph_\xi = \psi_\xi$.
Now $(\ref{needlet-rep})$ becomes
$f = \sum_{\xi \in \cX} \langle f, \psi_\xi\rangle \psi_\xi$.
It is easy to see that this representation holds
in $\L2$ and
$$ 
\|f\|_{\L2} =
\Big(\sum_{\xi \in \cX}
|\langle f, \psi_\xi\rangle|^2\Big)^{1/2},
\quad f\in \L2,
$$ 
i.e.
$\{\psi_\xi\}_{\xi \in \cX}$ is a tight frame for $\L2(\RR^d)$.
\end{re}

%
\section{Hermite-Triebel-Lizorkin spaces (F-spaces)}\label{Tri-Liz-spaces}
\setcounter{equation}{0}

In this section we introduce the analogue of Triebel-Lizorkin spaces
in the context of Hermite expansions following the general approach
described in \cite[\S 10.3]{T1}
and show that they can be characterized via needlets.
In our treatment of Hermite-Triebel-Lizorkin spaces we will utilize
the scheme of Frazier and Jawerth from \cite{F-J2} (see also \cite{F-J-W}).

\subsection{Definition of Hermite-Triebel-Lizorkin spaces}

Let the kernels $\{\Phi_j\}$ be defined by
\begin{equation}\label{def-Phi-j}
\Phi_0 := \CH_0
\quad\mbox{and}\quad
\Phi_j := \sum_{\nu=0}^\infty \ha
\Big(\frac{\nu}{4^{j-1}}\Big)\PP_\nu, \quad j\ge 1,
\end{equation}
where $\{\PP_\nu\}$ are from (\ref{def-Kn}) and
$\ha$ obeys the conditions:
\begin{align}
&\quad \ha\in C^\infty[0, \infty),
\quad \supp \, \ha \subset [1/4, 4], \label{ha1}\\
&\quad  |\ha(t)|>c>0, \quad \text{if } t \in [1/3, 3].\label{ha2}
\end{align}


\begin{defn}\label{def.Tr-Liz}
The Hermite-Triebel-Lizorkin space
$F^{\alpha q}_p:=F^{\alpha q}_p(H)$,
where $\alpha \in \R$, $0< p < \infty$, $0< q \le \infty$,
is defined as the set of all $f \in \cS'$ such that

\begin{equation}\label{F-norm}
\|f\|_{F^{\alpha q}_p} :=
\Big\|\Big(\sum_{j=0}^\infty
(2^{\alpha j}|\Phi_j*f(\cdot)|)^q\Big)^{1/q}\Big\|_p
< \infty,
\end{equation}
where the $\ell^q$-norm is replaced by the sup norm when $q=\infty$.
\end{defn}

As will be shown in
Theorem \ref{thm:needlet-Tr-Liz}, the above definition of
Triebel-Lizorkin spaces is independent of the specific selection
of $\ha$ satisfying (\ref{ha1})-(\ref{ha2}) in the definition
of $\Phi_j$ in (\ref{def-Phi-j}).


\begin{prop}\label{prop:ebedding}
The Hermite-Triebel-Lizorkin space $F^{\alpha q}_p$ is a quasi-Banach space
which is continuously embedded in $\cS'$
$($$F^{\alpha q}_p \hookrightarrow \cS'$$)$.
\end{prop}

\noindent
{\bf Proof.}
We will only establish that $F^{\alpha q}_p \hookrightarrow \cS'$.
Then the completeness of $F^{\alpha q}_p$ follows by a standard argument
using in addition Fatou's lemma and Proposition~\ref{prop:needlet-rep}.

As in Definition~\ref{def.Tr-Liz}, let $\{\Phi_j\}$ be
defined by a function $\ha$ obeying (\ref{ha1})-(\ref{ha2}).
As already indicated there exists a function $\hb$ such that
(\ref{ab1})-(\ref{ab5}) hold.
Let $\{\Psi_j\}$ be defined as in (\ref{def-Psi-j}) using this function.
After this preparation, let $\{\ph_\xi\}$ and $\{\psi_\xi\}$
be needlet systems defined as in (\ref{def-needlets1})-(\ref{def-needlets1})
using these $\{\Phi_j\}$ and $\{\Psi_j\}$.

Let $f\in F^{\alpha q}_p$. By Proposition~\ref{prop:needlet-rep}
$
f = \sum_{j=0}^\infty
\Psi_j*\overline{\Phi}_j*f
$
in $\cS'$ and hence
$$
\langle f, \phi \rangle = \sum_{j=0}^\infty
\langle \Psi_j*\overline{\Phi}_j*f, \phi \rangle
= \sum_{j=0}^\infty
\langle \overline{\Phi}_j*f, \overline{\Psi}_j*\phi \rangle,
\quad \phi\in\cS.
$$
Applying the Cauchy-Schwarz inequality and (\ref{norm-relation})
we obtain, for $j\ge 2$,
\begin{align*}
|\langle \overline{\Phi}_j*f, \overline{\Psi}_j*\phi \rangle|
&\le \|\Phi_j*f\|_2\|\Psi_j*\phi\|_2
\le c2^{jd/p}\|\Phi_j*f\|_p\sum_{\nu=4^{j-2}}^{4^j}\|\HH_j*\phi\|_2\\
&\le c2^{-j}\|f\|_{F^{\alpha q}_p}P_r^{*}(\phi),
\end{align*}
whenever $r\ge |\alpha|+d/p+1$.
This leads to
$
|\langle f, \phi\rangle|
\le c\|f\|_{F^{\alpha q}_p} P_r^{*}(\phi),
$
which yields the claimed embedding.
$\qed$


\begin{prop}\label{prop:identification}
We have the following identification:
\begin{equation}\label{ident1}
F^{0 2}_p \sim L^p,
\quad 1 < p < \infty,
\end{equation}
with equivalent norms.
\end{prop}
The proof of this proposition can be carried out as
the proof of Proposition~4.3 in \cite{NPW}
in the case of spherical harmonics and will be omitted.
It employs the existing $L^p$ multipliers for Hermite expansions
(see e.g. \cite{Th}).


\subsection{Needlet decomposition of Hermite-Triebel-Lizorkin spaces.}

In the following we will use the multilevel set $\cX:=\cup_{j=0}^\infty \cX_j$
from \S\ref{def-needlets} and the tiles $\{R_\xi\}$ introduced in (\ref{def.Q-xi}).


\begin{defn}\label{d:def.f-space}
Let $\alpha \in \R$, $0< p < \infty$, $0< q \le \infty$.
The Hermite-Triebel-Lizorkin sequence space
$f^{\alpha q}_p$ is defined as the set of all sequences of complex numbers
$s=\{s_\xi\}_{\xi \in \cX}$
such that
\begin{equation}\label{def-f-space}
\|s\|_{f^{\alpha q}_p} :=
\Big\|\Big(\sum_{j=0}^\infty 2^{j\alpha q}
\sum_{\xi \in \cX_j}
\Big[|s_\xi| |R_\xi|^{-1/2}
\ONE_{R_\xi}(\cdot)\Big]^q\Big)^{1/q}\Big\|_p <\infty
\end{equation}
with the usual modification when $q=\infty$.
\end{defn}

Assuming that $\{\ph_\xi\}$, $\{\psi_\xi\}$ is a dual pair
of analysis and synthesis needlets
(see (\ref{def-needlets1})-(\ref{def-needlets2})),
we introduce the operators:
$S_\ph: f \to \{\langle f, \ph_\xi\rangle\}_{\xi \in \cX}$
({\em Analysis operator})
and
$T_\psi: \{s_\xi\}_{\xi\in \cX} \to \sum_{\xi\in \cX}s_\xi\psi_\xi$
({\em Synthesis operator}).

We now come to our main result on Hermite-Triebel-Lizorkin spaces.


\begin{thm}\label{thm:needlet-Tr-Liz}
If $\alpha \in \R$ and $0<p<\infty$, $0<q\le\infty$,
then the operators
$S_\ph: F_p^{\alpha q} \to f_p^{\alpha q}$ and
$T_\psi: f_p^{\alpha q} \to F_p^{\alpha q}$ are bounded and
$T_\ph\circ S_\ph= {\rm Id}$.
Consequently, assuming that $f \in \cS'$, we have $f \in F_p^{\alpha q}$
if and only if $\{\langle f, \ph_\xi\rangle\}\in f_p^{\alpha q}$ and
\begin{equation}\label{Fnorm-needlet1}
\|f\|_{F_p^{\alpha q}}
\sim \|\{\langle f, \ph_\xi\rangle\}\|_{f_p^{\alpha q}}.
\end{equation}
Furthermore, the definition of $F_p^{\alpha q}$ is independent of the specific
selection of $\ha$ satisfying $(\ref{ha1})-(\ref{ha2})$.
\end{thm}

For the proof of this theorem we adapt some techniques from \cite{F-J2}.


\begin{defn}\label{def-s*}
For any collection of complex numbers $\{a_\xi\}_{\xi \in \cX_j}$,
we define
\begin{equation}\label{def-s}
a^*_j(x) := \sum_{\eta \in \cX_j} \frac{|a_\eta|}{(1+ 2^{j}|\eta-x|)^\sigma}
\end{equation}
and
\begin{equation}\label{def-a*}
a_\xi^*:=a^*_j(\xi),
\quad \xi \in \cX_j,
\end{equation}
where $\sigma >d$ is sufficiently large and will be specified later on.
\end{defn}

We will need a couple of lemmas whose proofs are given in \S\ref{proofs}.


\begin{lem}\label{lem:sum<M}
Suppose $s>0$ and $\sigma > d\max\{2, 1/s\}$.
Let $\{b_\omega\}_{\omega \in \cX_j}$, $j \ge 0$, be a set of complex numbers.
Then  
\begin{equation}\label{s*<M}
b_j^*(x)
\le c\cM_s\Big(\sum_{\omega \in \cX_j}|b_\omega| \ONE_{R_\omega}\Big)(x),
\quad x\in\RR^d.
\end{equation}
Moreover, for $\xi \in \cX_j$,
\begin{equation}\label{sum<M}
b_\xi^* \ONE_{R_\xi}(x)
\le c\cM_s\Big(\sum_{\omega \in \cX_j}|b_\omega| \ONE_{R_\omega}\Big)(x),
\quad x\in\RR^d.
\end{equation}
Here the constants depend only on $d$, $\delta$, $\sigma$, and $s$.
\end{lem}


\begin{lem}\label{lem:a*=b*}
Let $g \in \V_{4^j}$ and denote
$$
M_\xi:=\sup_{x \in R_\xi} |g(x)|,
\quad \xi\in\cX_j,
\quad\mbox{and}\quad
m_\lambda:=\inf_{x \in R_\lambda} |g(x)|, \quad
\lambda \in \cX_{j+\ell}.
$$
Then there exists $\ell \ge 1$, depending only $d$, $\delta$, and $\sigma$, such that
for any $\xi \in \cX_j$
\begin{equation}\label{a*=b*}
M_\xi^* \le c m_\lambda^*
\quad \mbox{for all \; $\lambda\in\cX_{j+\ell}$, $R_\lambda\cap R_\xi\ne \emptyset$},
\end{equation}
and hence
\begin{equation}\label{M*<m*2}
M_\xi^*\ONE_{R_\xi}(x)
\le c \sum_{\lambda\in\cX_{j+\ell}, R_\lambda\cap R_\xi\ne \emptyset}
m_\lambda^*\ONE_{R_\lambda}(x),
\quad x\in \R^d,
\end{equation}
where $c>0$ depends only on $d$, $\delta$, and $\sigma$.
\end{lem}

\medskip

\noindent
{\bf Proof of Theorem \ref{thm:needlet-Tr-Liz}.}
Suppose $q<\infty$ (the case $q=\infty$ is easier) and
pick $s$, $\sigma$, and $k$ so that
$0 < s < \min\{p, q\}$ and $k \ge \sigma > d\max\{1, 1/s\}$.

Let $\{\Phi_j\}$ be from the definition
of Hermite-Triebel-Lizorkin spaces (see (\ref{def-Phi-j})-(\ref{ha2})).
As already indicated in the beginning of \S\ref{def-needlets},
there exists a function $\hb$ satisfying (\ref{ab1})-(\ref{ab3})
such that (\ref{ab5}) holds as well.
We use this function to define $\{\Psi_j\}$ exactly as in (\ref{def-Psi-j}).
We further use $\{\Phi_j\}$ and $\{\Psi_j\}$ to define just as in (\ref{def-needlets1})
a~pair of dual needlet systems $\{\ph_\eta\}$ and $\{\psi_\eta\}$.

Let $\{\wtd\ph_\eta\}$, $\{\wtd\psi_\eta\}$ be a second pair of needlet systems,
defined as in (\ref{def.Phi-j})-(\ref{def-needlets1}) from another pair of
kernels $\{\wtd\Phi_j\}$, $\{\wtd\Psi_j\}$.

Our first step is to establish the boundedness of the operator
$T_{\wtd\psi}: f^{\alpha q}_p \to F^{\alpha q}_p$,
defined by
$
T_{\wtd\psi} s:= \sum_{\xi\in \cX} s_\xi\wtd\psi_\xi.
$
Proposition~\ref{prop:ebedding}
and the fact that finitely supported sequences are dense in $f^{\alpha q}_p$
imply that it suffices to prove the boundedness of $T_{\wtd\psi}$
only for finitely supported sequence.
So, assume $s=\{s_\xi\}_{\xi\in\cX}$ is a finitely supported sequence and
let $f:=T_{\wtd\psi}s$.
Evidently $\Phi_j*\wtd\psi_\xi =0$ if $\xi\in\cX_\nu$
and $|j-\nu|\ge 2$, and hence
$$
\Phi_j*f = \sum_{\nu=j-1}^{j+1}\sum_{\xi \in \cX_\nu}
s_\xi \Phi_j*\wtd\psi_\xi
\qquad (\cX_{-1}:=\emptyset).
$$

Let
$\xi \in \cX_{\nu}$, $j-1\le \nu \le j+1$, and $|\xi|_\infty \le (1+\delta)2^{\nu+1}$.
Then using (\ref{local-Needlets2})-(\ref{local-Needlets21}) we get
\begin{align*}
|\Phi_j*\wtd\psi_\xi(x)|
\le c2^{3jd/2}
\int_{\RR^d}\frac{1}{(1+2^j|x-y|)^k
(1+2^j|\xi-y|)^k}\, dy
\le \frac{c2^{jd/2}}{(1+2^j|\xi-x|)^k}.
\end{align*}
Hence, on account of (\ref{tiles1})
\begin{equation}\label{est-Phi-psi}
|\Phi_j*\wtd\psi_\xi(x)|
\le \frac{c|R_\xi|^{-1/2}}{(1+2^j|\xi-x|)^k},
\quad x\in \R^d.
\end{equation}

If $\xi \in \cX_{\nu}$, $j-1\le \nu \le j+1$, and $|\xi|_\infty > (1+\delta)2^{\nu+1}$,
then by (\ref{local-Needlets2})-(\ref{local-Needlets22})
\begin{align*}
|\Phi_j*\wtd\psi_\xi(x)|
\le c2^{-jL}
\int_{\RR^d}\frac{2^{jd}}{(1+2^j|x-y|)^k
(1+2^j|\xi-y|)^k}\, dy
\le \frac{c2^{-jL}}{(1+2^j|\xi-x|)^k}
\end{align*}
for any $L>0$.
Consequently, in view of (\ref{tiles2}), estimate (\ref{est-Phi-psi}) holds again.

Denote $S_\xi:= s_\xi|R_\xi|^{-1/2}$.
Then by (\ref{est-Phi-psi}) we have
\begin{align}\label{Phij-f<s*}
|\Phi_j*f(x)|
&\le \sum_{\nu=j-1}^{j+1}\sum_{\xi \in \cX_\nu}|s_\xi||\Phi_j*\wtd\psi_\xi(x)|
\le c\sum_{\nu=j-1}^{j+1}\sum_{\xi \in \cX_\nu}
\frac{|s_\xi||R_\xi|^{-1/2}}{(1+2^\nu |\xi-x|)^k}\notag\\
&\le c\sum_{\nu=j-1}^{j+1}S_\nu^*(x)
\qquad (S_{-1}:=0),
\end{align}
where
$S_\nu^*(x)$ is defined as in (\ref{def-s}).
We insert this in (\ref{F-norm}) and apply
Lemma \ref{lem:sum<M} and the maximal inequality (\ref{max-ineq}) to obtain
\begin{align*}
\|f\|_{F^{\alpha q}_p}
&\le \Big\|\Big(\sum_{j=0}^\infty
(2^{j\alpha}|S_j^*(\cdot)|)^q\Big)^{1/q}\Big\|_p\notag\\
&\le \Big\|\Big(\sum_{j=0}^\infty
\Big[\cM_s\Big(2^{j\alpha}\sum_{\xi\in\cX_j}|s_\xi||R_\xi|^{-1/2}\ONE_{R_\xi}\Big)
\Big]^q\Big)^{1/q}\Big\|_p
\le c\|\{s_\eta\}\|_{f^{\alpha q}_p}.
\end{align*}
Hence the operator $T_{\wtd\psi}: f^{\alpha q}_p \to F^{\alpha q}_p$
is bounded.


Assuming that the space $F^{\alpha q}_p$ is defined via $\{\overline{\Phi}_j\}$
instead of $\{\Phi_j\}$
we next prove the boundedness of the operator
$S_\ph: F^{\alpha q}_p \to f^{\alpha q}_p$.
Let $f\in F^{\alpha q}_p$ and set
$$
M_\xi := \sup_{x \in R_\xi}|\overline{\Phi}_j*f(x)|,
\quad \xi \in \cX_j,
\quad\mbox{and}\quad
m_\lambda :=\inf_{x \in R_\lambda}|\overline{\Phi}_j*f(x)|,
\quad \lambda \in \cX_{j+\ell},
$$
where $\ell$ is the constant from Lemma~\ref{lem:a*=b*}.
We have
\begin{eqnarray*}
|\langle f, \ph_\xi\rangle|
\le c|R_\xi|^{1/2}|\overline{\Phi}_j*f(\xi)|
\le c|R_\xi|^{1/2}M_\xi
\le c|R_\xi|^{1/2} M_\xi^*.
\end{eqnarray*}
By Lemma~\ref{lem:convolution}, $\overline{\Phi}_j*f \in \V_{4^j}$,
and applying Lemma \ref{lem:a*=b*} (see (\ref{M*<m*2})), we have
$$
M_\xi^*\ONE_{R_\xi}(x)
\le c \sum_{\lambda\in\cX_{j+\ell}, R_\lambda\cap R_\xi\ne \emptyset}
m_\lambda^*\ONE_{R_\lambda}(x),
\quad x\in \R^d.
$$
We use the above,
Lemma \ref{lem:sum<M}, and the maximal inequality (\ref{max-ineq})
to obtain
\begin{eqnarray*}
\|\{\langle f, \ph_\xi \rangle\}\|_{f^{\alpha q}_p}
&\le&
c\Big\|\Big(\sum_{j=0}^\infty 2^{\alpha jq}
\Big(\sum_{\xi \in \cX_{j}}
                 M_\xi^*\ONE_{R_\xi}\Big)^q\Big)^{1/q}\Big\|_p\\
&\le&
c\Big\|\Big(\sum_{j=0}^\infty 2^{\alpha jq}
\Big(\sum_{\lambda \in \cX_{j+\ell}}
                 m_\lambda^*\ONE_{R_\lambda}\Big)^q\Big)^{1/q}\Big\|_p\\
&\le& c\Big\|\Big(\sum_{j=0}^\infty
\cM_s\Big(2^{\alpha j}\sum_{\lambda \in \cX_{j+\ell}}
                   m_\lambda\ONE_{R_\lambda}\Big)^q\Big)^{1/q}\Big\|_p\\
&\le& c\Big\|\Big(\sum_{j=0}^\infty
\Big(2^{\alpha j}\sum_{\lambda \in \cX_{j+\ell}}
                   m_\lambda\ONE_{R_\xi}\Big)^q\Big)^{1/q}\Big\|_p\\
&\le& c\Big\|\Big(\sum_{j=0}^\infty 2^{\alpha jq}
|\overline{\Phi}_j*f|^q \Big)^{1/q}\Big\|_p
= c\|f\|_{F^{\alpha q}_p}.
\end{eqnarray*}
Here for the second inequality we used that each tile $R_\lambda$, $\lambda\in\cX_{j+l}$,
intersects no more that finitely many (depending only on $d$) tiles $R_\eta$, $\eta\in\cX_j$.
The above confirms the boundedness of the operator
$S_\ph: F^{\alpha q}_p \to f^{\alpha q}_p$.

The identity $T_\psi\circ S_\ph=Id$ follows by Theorem~\ref{prop:needlet-rep}.

We finally show the independence of the definition of Triebel-Lizorkin spaces
from the specific selection of $\ha$ satisfying (\ref{ha1})-(\ref{ha2}).
Let $\{\Phi_j\}$, $\{\wtd\Phi_j\}$ be two sequences of
kernels as in the definition of Triebel-Lizorkin spaces defined by
two different functions $\ha$ satisfying (\ref{ha1})-(\ref{ha2}).
As in the beginning of this proof, there exist two associated needlet systems
$\{\Phi_j\}$, $\{\Psi_j\}$, $\{\ph_\xi\}$, $\{\psi_\xi\}$
and
$\{\wtd\Phi_j\}$, $\{\wtd\Psi_j\}$, $\{\wtd\ph_\xi\}$, $\{\wtd\psi_\xi\}$.
Denote by $\|f\|_{F^{\alpha q}_p(\Phi)}$ and
$\|f\|_{F^{\alpha q}_p(\wtd\Phi)}$ the $F$-norms defined via
$\{\Phi_j\}$ and $\{\wtd\Phi_j\}$.
Then from above it follows that
$$
\|f\|_{F^{\alpha q}_p(\Phi)}
\le c \|\{\langle f, \wtd\ph_\xi\rangle\}\|_{f^{\alpha q}_p}
\le c\|f\|_{F^{\alpha q}_p(\overline{\wtd\Phi})}.
$$
The claimed independence of the definition of $F^{\alpha q}_p$ of the specific
selection of $\ha$ in the definition of the functions $\{\Phi_j\}$
follows by interchanging the roles of $\{\Phi_j\}$ and $\{\wtd\Phi_j\}$
and their complex conjugates.
\qed

\medskip

The Hermite-F-spaces embed in one another similarly as the classical F-spaces.

%
\begin{prop}\label{prop:F-embed-F}
$(a)$ If $0<p<\infty$, $0<q, q_1\le\infty$, $\alpha\in \R$ and
$\eps>0$, then
\begin{equation}\label{F-embedding1}
F^{\a+\eps, q}_{p} \hookrightarrow F^{\a q_1}_{p}.
\end{equation}

$(b)$ Let $0<p<p_1<\infty$, $0<q, q_1\le\infty$, and
$-\infty<\alpha_1<\alpha<\infty$. Then we have the continuous
embedding
\begin{equation}\label{F-embedding2}
\Fapq \hookrightarrow F^{\a_1q_1}_{p_1} \quad\mbox{if}\quad
\a-d/p=\a_1-d/p_1.
\end{equation}
\end{prop}
The proof of this embedding result uses estimate (\ref{norm-relation})
and Theorem~\ref{thm:needlet-Tr-Liz} and can
be carried out exactly as in the classical case on $\RR^n$
(see e.g. \cite{T1}, p. 47 and p.~129). We omit it.

\subsection{Comparison of Hermite-F-spaces with classical F-spaces}

We next use needlet decompositions to show that the
Hemite-Triebel-Lizorkin spaces of essentially
positive smoothness are different from the corresponding classical
Triebel-Lizorkin spaces on $\RR^d$.


\begin{thm}\label{thm:FH-diff-F}
Let $0<p <\infty$, $0<q\le\infty$, and $\alpha>d(1/p-1)_+$.
Then there exists a function $f\in F^{\a q}_p$
such that $\|f\|_{F_p^{\a q}(H)}=\infty$ and hence
$f\not\in F^{\a q}_p(H)$.
Here $F^{\a q}_p$  and $F^{\a q}_p(H)$ are the respective
classical and Hermite Triebel-Lizorkin spaces.
\end{thm}

\noindent
{\bf Proof.}
For any $y\in\R^d$ and a function $f$ we define
\begin{equation}\label{def-F*}
\|f\|_{F^*_y}:=
\Big\|\Big(  \sum_{j=0}^\infty  2^{j\a q}\sum_{\xi\in\cX_j,\, |\xi-y|>|y|/2}
\Big(|R_\xi|^{-1/2}|\langle f, \ph_\xi\rangle|
\ONE_{R_\xi}(\cdot)\Big)^q\Big)^{1/q}\Big\|_p.
\end{equation}
Choose a function $h\in C^\infty(\R^d)$ such that $\|h\|_\infty=1$ and
$\supp h \subset B(0, 1)$, where $B(0, 1):=\{x\in\R^d: |x|<1\}$.

Theorem~\ref{thm:FH-diff-F} will follow easily by the following lemma
whose proof is given in \S\ref{proofs1}.


\begin{lem}\label{lem:FH-diff-F}
With the notation from above, we have
\begin{align}
\|h(\cdot - y)\|_{F^{\a q}_p(H)} \to \infty
\quad &\mbox{ as } \quad
|y| \to \infty, \quad \mbox{and}\label{FH-diff-F1}\\
\|h(\cdot - y)\|_{F^*_y} \to 0
\quad &\mbox{ as } \quad
|y| \to \infty. \label{FH-diff-F2}
\end{align}
\end{lem}


By this lemma it follows that there exists a sequence
$\{y_j\}_{j\ge 1}\subset \R^d$ such that
$0<|y_1|<|y_2| < \dots$ and $|y_{j+1}|>3|y_j|$,
$\|h(\cdot-y_j)\|_{F_q^{\a p}(H)}>2^{2j}$,
and
$\|h(\cdot-y_j)\|_{F^*_{y_j}}<1$,
$j=1, 2, \dots$.

We now define
$f(x):= \sum_{j=1}^\infty f_j(x)$, where $f_j(x):=2^{-j}h(x-y_j)$, $x\in\R^d$.
Set $\tau:=\min\{p, q, 1\}$.
Evidently, $h$ belongs to all classical Triebel-Lizorkin spaces,
which are shift invariant, and hence
$$
\|f\|_{F_q^{\a p}}^\tau
\le \sum_{j=1}^\infty 2^{-j\tau}\|h(\cdot-y_j)\|_{F_q^{\a p}}^\tau
=\|h\|_{F_q^{\a p}}^\tau \sum_{j=1}^\infty 2^{-j\tau}
\le c\|h\|_{F_q^{\a p}}^\tau<\infty.
$$
Here we use that
$\|\sum_j g_j\|_{F_q^{\a p}}^\tau \le \sum_j\|g_j\|_{F_q^{\a p}}^\tau.$
Thus $f\in F_q^{\a p}$.

On the other hand, for any $\ell\ge 1$,
\begin{align*}
\|f\|_{F_q^{\a p}(H)}^\tau
&\ge c\Big\|\Big(  \sum_{\nu=0}^\infty  2^{\nu\a q}\sum_{\xi\in\cX_\nu,\, |\xi-y_\ell|\le|y_\ell|/2}
\Big(|R_\xi|^{-1/2}|\langle f, \ph_\xi\rangle|
\ONE_{R_\xi}(\cdot)\Big)^q\Big)^{1/q}\Big\|_p^\tau\\
&\ge c\Big(\|f_\ell\|_{F_q^{\a p}(H)}^\tau
-\sum_{j=1}^\infty\|f_j\|_{F^*_{y_j}}^\tau\Big)\\
&= c2^{-\ell\tau}\|h(\cdot-y_\ell)\|_{F_q^{\a p}(H)}^\tau
-c\sum_{j=1}^\infty 2^{-j\tau}\|h(\cdot-y_j)\|_{F^*_{y_j}}^\tau\\
&> c2^{\ell\tau}- c\sum_{j=1}^\infty 2^{-j\tau} \ge c2^{\ell\tau}-c'.
\end{align*}
Here for the second inequality we used that if $|\xi-y_\ell|\le |y_\ell|/2$,
then $|\xi-y_j|>|y_j|/2$ for all $j\ne \ell$.
Consequently, $\|f\|_{F_q^{\a p}(H)}^\tau =\infty$.
$\qed$



\section{Hermite-Besov spaces (B-spaces)}\label{Besov-spaces}
\setcounter{equation}{0}

Besov type spaces are natural to introduce in the context of Hermite expansions
(see e.g. \cite[\S 10.3]{T1}).
We will call them Hermite-Besov spaces.
To characterize these space via needlets we use the approach of
Frazier and Jawerth \cite{F-J1} (see also \cite{F-J-W})
to the classical Besov spaces.
We refer to \cite{Peetre, T1} as general references for Besov spaces.

\subsection{Definition of Hermite-Besov spaces}


\begin{defn}\label{def.Besov-sp}
Let the kernels $\{\Phi_j\}$ be defined by $(\ref{def-Phi-j})$ with $\ha$
satisfying $(\ref{ha1})-(\ref{ha2})$.
The Hermite-Besov space $B^{\alpha q}_p := B^{\alpha q}_p(H)$,
where $\alpha \in \R$, $0<p,q \le \infty$, is defined
as the set of all $f \in \cS'$ such that
\begin{equation}\label{B-norm}
\|f\|_{B^{\alpha q}_p} :=
\Big(\sum_{j=0}^\infty \Big(2^{\alpha j}\|\Phi_j*f\|_p\Big)^q\Big)^{1/q}
< \infty,
\end{equation}
where the $\ell^q$-norm is replaced by the sup-norm if $q=\infty$.
\end{defn}

Similarly as for Hermite-Triebel-Lizorkin spaces (\S\ref{Tri-Liz-spaces})
Theorem \ref{thm:needlet-Besov} below implies that the above definition
of Hermite-Besov spaces is independent of the specific selection of $\ha$;
also $B^{\alpha q}_p$ is a quasi-Banach space which is continuously
embedded in $\cS'$.


\subsection{Needlet decomposition of Hermite-Besov spaces}

As for the Hermite-Triebel-Lizorkin spaces we employ
the tiles $\{R_\xi\}$ introduced in (\ref{def.Q-xi}) in the following.
Also as before $\cX:=\cup_{j=0}^\infty \cX_j$.


\begin{defn}\label{def.b-space}
The Hermite-Besov sequence space $\bb^{\alpha q}_p$,
where $\alpha \in \R$, $0<p,q \le \infty$,
is defined
as the set of all sequences of complex numbers
$s=\{s_\xi\}_{\xi \in \cX}$
such that
\begin{equation}\label{def-b-space}
\|s\|_{\bb^{\alpha q}_p} :=
\Big(\sum_{j=0}^\infty \Big[2^{j\alpha}\Big(\sum_{\xi \in \cX_j}
|R_\xi|^{1-p/2}|s_\xi|^p\Big)^{1/p}\Big]^{q}\Big)^{1/q} <\infty
\end{equation}
with obvious modifications when $p=\infty$ or $q=\infty$.
\end{defn}

In the following, we assume that
$\{\Phi_j\}$, $\{\Psi_j\}$, $\{\ph_\xi\}$, $\{\psi_\xi\}$
is a needlet system defined by (\ref{def.Phi-j})-(\ref{def-needlets2}).
Recall the
{\em analysis operator:}
$S_\ph: f \to \{\langle f, \ph_\xi\rangle\}_{\xi \in \cX}$,
and the
{\em synthesis operator:}
$T_\psi: \{s_\xi\}_{\xi\in \cX} \to \sum_{\xi\in \cX}s_\xi\psi_\xi$.


\begin{thm}\label{thm:needlet-Besov}
If $\alpha \in \R$ and $0<p, q\le\infty$,
then the operators
$S_\ph: B_p^{\alpha q} \to \bb_p^{\alpha q}$ and
$T_\psi: \bb_p^{\alpha q} \to B_p^{\alpha q}$ are bounded and
$T_\ph\circ S_\ph= {\rm Id}$.
Consequently,
assuming that $f \in \cS'$, we have $f \in B_p^{\alpha q}$
if and only if $\{\langle f, \ph_\xi\rangle\}\in \bb_p^{\alpha q}$ and
\begin{equation}\label{Bnorm-needlet1}
\|f\|_{B_p^{\alpha q}}
\sim \|\{\langle f, \ph_\xi\rangle\}\|_{\bb_p^{\alpha q}}.
\end{equation}
Furthermore, the definition of $B_p^{\alpha q}$ is independent of
the choice of $\ha$ satisfying $(\ref{ha1})-(\ref{ha2})$.
\end{thm}

For the proof of this theorem we need one additional lemma.


\begin{lem}\label{lem:Polyn-Lp}
For any $g\in \V_{4^j}$, $j\ge 0$, and $0<p\le \infty$
\begin{equation}
 \Big(\sum_{\xi\in \cX_j}|R_\xi|\max_{x\in R_\xi}|g(x)|^p\Big)^{1/p}
 \le c\|g\|_p.
\end{equation}
\end{lem}

The proof of this lemma is given in \S\ref{proofs}.

\medskip


\noindent
{\bf Proof of Theorem~\ref{thm:needlet-Besov}.}
Let $0<s<p$ and $\sigma >d\max\{1, 1/s\}$.
Just as in the proof of Theorem~\ref{thm:needlet-Tr-Liz}
we assume that
$\{\Phi_j\}$, $\{\Psi_j\}$, $\{\ph_\eta\}$, $\{\psi_\eta\}$
and
$\{\wtd\Phi_j\}$, $\{\wtd\Psi_j\}$, $\{\wtd\ph_\eta\}$, $\{\wtd\psi_\eta\}$
are two needlet systems, defined as in (\ref{def.Phi-j})-(\ref{def-needlets1}),
which originate from two completely different functions $\ha$ satisfying
(\ref{ha1})-(\ref{ha2}).

We first prove the boundedness of the operator
$T_{\wtd\psi}: \bapq \to \Bapq$,
defined by
$
T_{\wtd\psi} s:= \sum_{\xi\in \cX} s_\xi\wtd\psi_\xi,
$
assuming that $\Bapq$ is defined by $\{\Phi_j\}$.
As in the Triebel-Lizorkin case due to the embedding $\Bapq \hookrightarrow \cS'$
it suffices to consider only the case of a finitely supported sequence
$s=\{s_\xi\}_{\xi\in\cX}$. Let $f:=T_{\wtd\psi}s.$
By (\ref{Phij-f<s*})
and Lemma~\ref{lem:sum<M} we get
\begin{align*}
\|\Phi_j*f\|_p
&\le c\sum_{\nu=j-1}^{j+1}
\Big\|\cM_s\Big(\sum_{\omega\in\cX_\nu}|R_\omega|^{-1/2}|s_\omega|\ONE_{R_\omega}\Big)\Big\|_p\\
&\le c\sum_{\nu=j-1}^{j+1}
\Big(\sum_{\omega\in\cX_\nu}|R_\omega|^{1-p/2}|s_\omega|^p\Big)^{1/p}
\quad (\cX_{-1}:=\emptyset),
\end{align*}
which leads to
$\|f\|_{\Bapq} \le c \|\{s_\eta\}\|_{\bapq}$
and hence to the boundedness of $T_{\wtd\psi}$.

To prove the boundedness of the operator $S_\varphi:\Bapq\rightarrow \bapq$
we assume that $\Bapq$ is defined in terms of $\{\overline{\Phi}_j\}$.
Observing that
$$
|\langle f,\varphi_\xi \rangle|
=\lambda_\xi^{1/2}|\overline{\Phi}_j*f(\xi)|
\sim |R_\xi|^{1/2}|\overline{\Phi}_j*f(\xi)|,
\quad \xi \in \cX_j,
$$
and $\overline{\Phi}_j*f\in \V_{4^j}$, we get using Lemma~\ref{lem:Polyn-Lp}
\begin{align*}
\Big(\sum_{\xi\in\cX_j}|R_\xi|^{1-p/2}|\langle f,\varphi_\xi \rangle|^p\Big)^{1/p}
\le c\Big(\sum_{\xi\in\cX_j}|R_\xi||\overline{\Phi}_j*f(\xi)|^p\Big)^{1/p}
\le c\|\overline{\Phi}_j*f\|_p.
\end{align*}
This yields
$
\|\{\langle f,\varphi_\xi \rangle\}\|_{\bapq} \le c \|f\|_{\Bapq}
$
and hence the operator $S_\varphi$ is bounded.

The identity $T_\psi\circ S_\varphi=Id$ is a consequence of
 Proposition~\ref{prop:needlet-rep}.

The independence of the definition of $\Bapq$ from the particular selection
of $\ha$ follows from above exactly as in the case of Triebel-lizorkin spaces
(see the proof of Theorem~\ref{thm:needlet-Tr-Liz}).
$\qed$

\medskip

The Hermite-Besov spaces embed similarly as the classical Besov spaces.

%
\begin{prop}\label{prop:B-embed-B}
$(a)$ If $0<p, q, q_1\le\infty$, $\alpha\in \R$ and
$\eps>0$, then
\begin{equation}\label{B-embedding1}
B^{\a+\eps, q}_{p} \hookrightarrow B^{\a q_1}_{p}.
\end{equation}

$(b)$ Let $0<p < p_1<\infty$, $0<q\le\infty$, and
$-\infty<\alpha_1<\alpha<\infty$. Then we have the continuous
embedding
\begin{equation}\label{B-embedding2}
\Bapq \hookrightarrow B^{\a_1q}_{p_1} \quad\mbox{if}\quad
\a-d/p=\a_1-d/p_1.
\end{equation}

$(c)$ If $0<p<\infty$, $0<q\le\infty$, $\alpha\in \R$, then
\begin{equation}\label{B-F-embedding}
B^{\a, \min\{p, q\}}_{p} \hookrightarrow F^{\a q}_{p}
\hookrightarrow B^{\a, \max\{p, q\}}_{p}.
\end{equation}

\end{prop}

Part (b) of this proposition follows readily by estimate (\ref{norm-relation}).
The proofs of parts (a) and (c) are as in the classical case.

\medskip

We now show that under some restriction on the indices the
Hemite-Besov spaces are essentially different from the classical
Besov spaces on $\RR^d$.


\begin{thm}\label{thm:B-ne-B}
Let $0<p, q <\infty$, and $\alpha>d(1/p-1)_+$.
Then there exists a function $f\in B^{\a q}_p$
such that $\|f\|_{B^{\a q}_p(H)}=\infty$ and then
$f\not\in B^{\a q}_p(H)$.
Here $B^{\a q}_p$  and $B^{\a q}_p(H)$ are the respective
classical and Hermite Besov spaces.
\end{thm}

\noindent
{\bf Proof.}
We proceed quite similarly as in the proof of Theorem~\ref{thm:FH-diff-F}.
Given $y\in\R^d$ and a function $f$ we define
\begin{equation}\label{def-B*}
\|f\|_{B^*_y}:=
\Big(  \sum_{j=0}^\infty  2^{j\a q}
\Big(\sum_{\xi\in\cX_j,\, |\xi-y|>|y|/2}
|R_\xi|^{1-p/2}|\langle f, \ph_\xi\rangle|^p
\Big)^{q/p}\Big)^{1/q}.
\end{equation}
Pick $h\in C^\infty(\R^d)$ such that $\|h\|_\infty=1$ and
$\supp h \subset B(0, 1)$.

The theorem follows easily by the following:
\begin{align}
\|h(\cdot - y)\|_{B^{\a q}_p(H)} \to \infty
\quad &\mbox{ as } \quad
|y| \to \infty, \quad \mbox{and}\label{B-ne-B1}\\
\|h(\cdot - y)\|_{B^*_y} \to 0
\quad &\mbox{ as } \quad
|y| \to \infty. \label{B-ne-B2}
\end{align}

To prove (\ref{B-ne-B1}) we will show that there exist $\eps>0$ and
$r>1$ such that
\begin{equation}\label{embed-B-F}
B^{\a q}_p(H) \hookrightarrow F^{\eps 2}_{r}(H).
\end{equation}
Then the result follows by the argument from the proof of
Lemma~\ref{lem:FH-diff-F}.

Let $p>1$. Pick $\eps>0$ so that $\alpha > 2\eps$.
Then by Propositions~\ref{prop:F-embed-F}, \ref{prop:B-embed-B}
we have the following embeddings
$
B^{\a q}_p \hookrightarrow B^{2\eps, p}_p
           \hookrightarrow F^{2\eps, p}_p
           \hookrightarrow F^{\eps 2}_p
$
which confirms (\ref{embed-B-F}).

Let $p\le 1$. Then $\alpha > d(1-1/p)$ and hence
(as in the proof of Lemma~\ref{lem:FH-diff-F}) there exist
$\eps, \delta>0$ such that
$\alpha-d/p=3\eps-d/(1+\delta)$.
Then Propositions~\ref{prop:F-embed-F}, \ref{prop:B-embed-B}
give us the following embeddings
$
B^{\a q}_p \hookrightarrow B^{3\eps, q}_{1+\delta}
           \hookrightarrow B^{2\eps, 1+\delta}_{1+\delta}
           \hookrightarrow F^{2\eps, 1+\delta}_{1+\delta}
           \hookrightarrow F^{\eps 2}_{1+\delta},
$
which leads again to (\ref{embed-B-F}).

The proof of (\ref{B-ne-B2}) is similar to the proof
of (\ref{FH-diff-F2}) and will be omitted.
$\qed$

\medskip

%
We finally want to link the Hermite-Besov spaces
with the $L^p$-approximation from linear combinations of Hermite functions.
Denote by $E_n(f)_p$ the best approximation of $f \in L^p$ from
$\V_n$, i.e.
\begin{equation}\label{def-Em}
E_n(f)_p := \inf_{g \in \V_n}\|f-g\|_p.
\end{equation}
Let $A_p^{\alpha q}$ be the approximation space of all functions $f\in L^p$ for which
\begin{equation}\label{App-space}
\|f\|_{A_p^{\alpha q}}
:= \|f\|_p +
\Big(\sum_{j=0}^\infty (2^{\alpha j}E_{2^j}(f)_p)^q\Big)^{1/q}
< \infty
\end{equation}
with the usual modification when $q=\infty$.


\begin{prop}\label{thm:character-Besov}
If $\alpha > 0$, $1 \le p \le \infty$, $0 < q \le \infty$,
then $B^{\alpha q}_p=A_p^{\alpha/2, q}$
with equivalent norms.
\end{prop}

\noindent
{\bf Proof.}
Let $f \in B^{\alpha q}_p$.
We first observe that under the conditions on $\alpha$, $p$, and $q$,
$B^{\alpha q}_p$ is continuously imbedded in $L^p$,
i.e. $f$ can be identified as a function in $L^p$ and
$\|f\|_p \le c\|f\|_{B^{\alpha q}_p}$.
The proof of this is easy and standard and will be omitted.

By a well known and easy construction there exists a function
$\ha\ge 0$ satisfying (\ref{ha1})-(\ref{ha2}) such that
$\ha(t)+\ha(4t)=1$ for $t\in [1/4, 1]$ and hence
$\sum_{\nu=0}^\infty \ha(4^{-\nu}t) =1$ for $t\in [1, \infty)$.
Assume that $\{\Phi_j\}$ are defined by (\ref{def-Phi-j})
using this function $\ha$.
By~Theorem~\ref{thm:needlet-Besov} the definition of the Besov spaces
$B^{\alpha q}_p$ is independent of the selection of $\ha$ and hence
they can be defined via these functions $\{\Phi_j\}$.
Similarly as in Proposition~\ref{prop:needlet-rep}
$f = \sum_{j=0}^\infty \Phi_j*f$ for $f\in L^p$.

Now, using that $\Phi_j*f \in \V_{4^{j}}$, we have
$
E_{4^m}(f)_p \le \sum_{j=m+1}^\infty \|\Phi_j*f\|_p,
$
$m\ge 0$.
A~standard argument employing this leads to the estimate
$\|f\|_{A^{\alpha/2, q}_p}^A \le c\|f\|_{B^{\alpha q}_p}$.

To prove the estimate in the other direction,
let $g \in \V_{4^{j-2}}$ ($j\ge 2$).
Evidently,
$\Phi_j*f = \Phi_j*(f-g)$
and the rapid decay of $\Phi_j$ yields
$
\|\Phi_j*f\|_p \le c\|f-g\|_p.
$
Consequently,
$\|\Phi_j*f\|_p \le cE_{4^{j-2}}(f)_p$, $j\ge 2$,
and
$\|\Phi_j*f\|_p \le c\|f\|_p$.
These lead to
$
\|f\|_{B^{\alpha q}_p} \le c\|f\|_{A^{\alpha/2, q}_p}^A.
$
$\qed$


\section{Proofs}\label{proofs}
\setcounter{equation}{0}

\subsection{Proofs for Section~\ref{Local-kernels}.}\label{proof-loc-est}
${}$

\medskip

\noindent
{\bf Proof of Theorem \ref{thm:derivative}.}
This proof hinges on an important lemma from \cite{Th}.
Let $\psi$ be a univariate function. The forward differences of
$\psi$ are defined by
$$
\Delta \psi(t) = \psi(t+1)-\psi(t)
\quad\mbox{and}\quad
\Delta^k \psi =
\Delta (\Delta^{k-1} \psi), \quad k \ge 2.
$$
For a given function $\psi$ we define
$$
M_\psi(x,y) := \sum_{\nu=0}^\infty \psi(\nu) \HH_\nu(x,y), \quad\hbox{and then}\quad
M_{\Delta^k \psi} = \sum_{\nu=0}^\infty \Delta^k\psi(\nu) \HH_\nu.
$$


\begin{lem}\label{lem:main}\cite[p. 72]{Th}
Let $A_j^{(x)}$ and $A_j^{(y)}$ denote the operator $A_j$ applied to the $x$ and
$y$ variables, respectively. Then for any $k \ge 1$,
\begin{equation}\label{eq:M-psi}
2^k (x_j-y_j)^k M_\psi(x,y) = \sum_{k/2 \le l \le k} c_{l,k}
\left(A_j^{(y)} - A_j^{(x)}\right)^{2l-k} M_{\Delta^l\psi}(x,y),
\end{equation}
where $c_{l,k}$ are constants given by
$$
c_{l,k} = (-1)^{k-l} 4^{k-l} (2k-2l-1)!! \binom{k}{2l-k}.
$$
\end{lem}

\noindent
{\bf Proof.}
This lemma is proved in \cite{Th} except that the constants $c_{l,k}$ are not
determined explicitly there. For $k=1$ one has \cite[(3.2.20)]{Th}
\begin{equation}\label{eq:k=1}
2 (x_j - y_j) M_\psi(x,y) = \left(A_j^{(y)} -A_j^{(x)}\right) M_{\Delta\psi}(x,y).
\end{equation}
The general result is obtained by induction using the identity \cite[(3.2.23)]{Th}
\begin{equation}\label{k-general}
(x_j - y_j) \left( A_j^{(y)} -A_j^{(x)}\right)^r -
\left( A_j^{(y)} -A_j^{(x)}\right)^r (x_j-y_j) = - 2 r
\left( A_j^{(y)} -A_j^{(x)}\right)^{r-1}.
\end{equation}
Assume that (\ref{eq:M-psi}) holds for some $k\ge0$.
Then using (\ref{k-general}) we get
\begin{align*}
& 2^{k+1} (x_j - y_j)^{k+1} M_\psi = 2(x_j-y_j)
\sum_{k/2 \le l \le k} c_{l,k}
\left(A_j^{(y)} - A_j^{(x)}\right)^{2l-k} M_{\Delta^l\psi} \\
& \qquad = \sum_{k/2 \le l \le k} c_{l,k} \left[
\left(A_j^{(y)} - A_j^{(x)}\right)^{2l-k}2(x_j-y_j) M_{\Delta^l\psi} \right.\\
& \qquad\qquad \left.- 4(2l-k)
\left(A_j^{(y)}-A_j^{(x)}\right)^{2l-k-1} M_{\Delta^l\psi} \right]\\
& \qquad= c_{k,k} \left(A_j^{(y)} - A_j^{(x)}\right)^{k+1}M_{\Delta^{k+1}\psi} \\
& \qquad\qquad + \sum_{(k+1)/2 \le l \le k} \left[c_{l-1,k} - 4 (2l-k)c_{l,k}\right]
\left(A_j^{(y)} - A_j^{(x)}\right)^{2l-k-1} M_{\Delta^l\psi},
\end{align*}
where $c_{l,k}:=0$ if $l<k/2$.
Consequently, the coefficients satisfy the recurrence relations
$$
c_{k+1,k+1}=c_{k,k}, \qquad c_{l,k+1}=c_{l-1,k} - 4 (2l-k)c_{l,k},
\quad (k+1)/2 \le l \le k.
$$
It follows from this and (\ref{eq:k=1}) that $c_{k,k} = 1$ for all $k$.
Furthermore, the recurrence relation above shows that
$$
c_{k-j,k} = -4 \sum_{\nu = 2j-1}^{k-1} (\nu - 2j+2)c_{\nu-j+1,\nu},
$$
from which one uses induction and the fact that
$
\sum_{\nu=j}^{k-1} \binom{\nu}{j} = \binom{k}{j+1}
$
to derive the stated identity for $c_{l,k}$.
$\qed$

\medskip


{\bf The case \boldmath $\alpha=(0, \dots, 0)$.}
Assume $k\ge 2$.
By~Lemma~\ref{lem:main}, we have
\begin{align*}
2^k(x_j-y_j)^k \LL_n(x,y) = \sum_{k/2 \le l \le k} c_{l,k}
\sum_{\nu=0}^{\infty} \Delta^l \ha\Big(\frac{\nu}{n}\Big)
\left(A_j^{(y)} - A_j^{(x)}\right)^{2l-k} \Ph_\nu(x,y),
\end{align*}
where $\Delta^l \ha(\frac{\nu}{n})$ is the $l$th forward difference applied with
respect to $\nu$.

Note first that applying the Cauchy-Schwarz inequality, we have
\begin{align}\label{H<HH}
\Big |\sum_{|\alpha|= n} \HH_{\alpha+\beta}(x) \HH_{\alpha+\gamma}(y) \Big|^2
    &\le \sum_{|\alpha|=n} |\HH_{\alpha+\beta}(x)|^2
         \sum_{|\alpha|=n} |\HH_{\alpha+\gamma}(y)|^2 \notag\\
  &  \le \HH_{n+|\beta|}(x,x) \HH_{n+|\gamma|}(y,y).
\end{align}

We next consider the action of $A_j$ on $\Ph_j(x,y)$.
Using repeatedly (\ref{eq:Aj}) we get
\begin{equation}\label{AjHH}
\Big(A_j^{(x)}\Big)^m\HH_\beta
= \prod_{r=0}^m[2(\beta_j+r)+2]^{1/2} \HH_{\beta+me_j}.
\end{equation}
This along with (\ref{H<HH}) leads to
\begin{align*}
|(A_j^{(y)})^i (A_j^{(x)})^{m}\Ph_\nu(x,y)|^2
&=\Big| \sum_{|\alpha|=\nu}
\prod_{r=0}^{i-1}[2(\alpha_j+r)+2]^{1/2}\prod_{r=0}^{m-1}[2(\alpha_j+r)+2]^{1/2} \notag \\
&\qquad\qquad \times \Ph_{\a + i e_j}(y) \Ph_{\a + me_j}(x) \Big|^2
   \notag  \\
& \le [2(\nu+i+m)]^{i+m} \left[ \sum_{|\alpha|=\nu}
    \left | \Ph_{\a + i e_j}(y) \Ph_{\a + me_j}(x) \right| \right]^2 \\
& \le [2(\nu+i+m)]^{i+m}  \Ph_{\nu + i}(y,y) \Ph_{\nu + m}(x,x).\notag
\end{align*}
Hence the binomial theorem and the Cauchy-Schwarz inequality give
\begin{align*}
& \Big| \left(A_j^{(y)} - A_j^{(x)}\right)^{2l-k} \Ph_\nu(x,y)\Big|
\le \sum_{i=0}^{2l-k} \binom{2l-k}{i}
\Big|(A_j^{(y)})^i (A_j^{(x)})^{2l-k-i}\Ph_\nu(x,y)\Big|\\
& \qquad \le (2\nu+4l-2k)^{(2l-k)/2} \sum_{i=0}^{2l-k} \binom{2l-k}{i}
   \left[\Ph_{\nu + i}(y,y)\right]^{\frac12}
      \left[ \Ph_{\nu + 2l-k-i}(x,x)\right]^{\frac12} \\
& \qquad \le c \nu^{(2l-k)/2}
      \left[\sum_{i=0}^{2l-k} \Ph_{\nu + i}(y,y)\right]^{\frac12}
     \left[\sum_{i=0}^{2l-k} \Ph_{\nu + i}(x,x) \right]^{\frac12}.
\end{align*}
A well known property of the difference operator gives
\begin{equation} \label{eq:hat-a}
\Big|\Delta^l \wh a \Big(\frac{\nu}{n}\Big) \Big|
  = n ^{-l} |\wh a^{(l)} (\xi) | \le n^{-l}\|\wh a^{(l)}\|_{\Linfty}.
\end{equation}
By Definition~\ref{defn:admissible}
it follow that $\Delta^l \ha(\frac{\nu}{n})=0$ if
$0\le \nu\le un-l$ or $\nu\ge n+vn$, where $0<u\le 1$ and $v>0$.
(Here $u=1$ if $\ha$ is of type (a).)
Using this, the above estimates, and the Cauchy-Schwarz inequality we infer
\begin{align*}
\left|(x_j-y_j)^k \LL_n(x,y)\right| & \le c
\sum_{k/2 \le l \le k}|c_{l,k}| n^{-l} \|\wh a^{(l)}\|_\infty n^{(2l-k)/2} \\
& \quad\qquad  \times
    \sum_{ \nu=[un]-l}^{n+[vn]}  \left[\sum_{i=0}^{2l-k} \Ph_{\nu + i}(y,y)\right]^{\frac12}
     \left[\sum_{i=0}^{2l-k} \Ph_{\nu + i}(x,x) \right]^{\frac12}  \\
%
%
& \le c_k n^{-k /2} \left[K_{n+[vn]+k}(y,y)\right]^{\frac12}
    \left[K_{n+[vn] + k}(x,x)\right]^{\frac12},
\end{align*}
where $c_k>0$ is of the form $c_k=c(k, u, d)\max_{0\le l\le k}\|\ha^{(l)}\|_\infty$.
Consequently,
\begin{equation} \label{mid-est2}
|\LL_n(x,y)|\le c_k\frac{\left[K_{n+[vn]+k}(y,y)\right]^{\frac12}
    \left[K_{n+[vn]+k}(x,x)\right]^{\frac12} }{(\sqrt{n}|x-y|)^{k}},
\qquad x\ne y.
\end{equation}

We also need estimate $|\LL_n(x,y)|$ whenever $x$ and $y$ are close to
one another. Applying the Cauchy-Schwarz inequality to the sum
in (\ref{def-Kn}) which defines $\Ph_\nu$ we get
\begin{align*}
|\LL_n(x,y)| &\le \sum_{\nu=0}^{n+[vn]}
\Big|\ha\Big(\frac{\nu}{n}\Big)\Big|
\Ph_\nu(x,x)^{1/2} \Ph_\nu(y,y)^{1/2} \\
& \le c \|\ha\|_\infty
  \left[K_{n+[vn]}(y,y)\right]^{\frac12} \left[K_{n+[vn]}(x,x)\right]^{\frac12},
\end{align*}
which coupled with (\ref{mid-est2}) yields (\ref{eq:derivative1})
in the case under consideration.

\medskip

{\bf The case \boldmath $|\alpha|>0$.}
We will make use of the relation $\partial_j = x_j - A_j$, where as usual
$\partial_j:=\frac{\partial}{\partial x_j}$.
In the following we
again denote by $A_j^{(x)}$ the operator $A_j$ acting on the $x$ variables, and
$A_j^0$ is understood as the identity operator;
by~definition
$(A^{(x)})^\alpha:=(A_1^{(x)})^{\alpha_1}\dots (A_d^{(x)})^{\alpha_d}$.
We also identify the operator of multiplication by $x_j$ with $x_j$. In order to use
Lemma \ref{lem:main}, we will need two commuting relations.


\begin{lem}\label{lem:comm}
Let $k, r, s$ be nonnegative integers. Then
\begin{equation} \label{eq:comm1}
    x_j^r \left(A_j^{(x)} - A_j^{(y)}\right)^k = \sum_{i=0}^r \binom{r}{i} \frac{k!}{(k-i)!}
       \left(A_j^{(x)} - A_j^{(y)}\right)^{k-i} x_j^{r-i}
\end{equation}
and
\begin{equation}\label{eq:comm2}
  (x_j-y_j)^k \left(A_j^{(x)}\right)^s = \sum_{i=0}^s \binom{s}{i} \frac{k!}{(k-i)!}
       \left(A_j^{(x)}\right)^{s-i} (x_j-y_j)^{k-i},
\end{equation}
where $k!/(k-i)! := 0$ if $k < i$.
\end{lem}

\noindent
{\bf Proof.}
To prove (\ref{eq:comm1}) we start from the identity:
\begin{equation} \label{xjA=Axj}
x_j  \left(A_j^{(x)}-  A_j^{(y)}\right)^k =  k  \left(A_j^{(x)}-  A_j^{(y)}\right) ^{k-1} +
\left( A_j^{(x)}- A_j^{(y)} \right)^k x_j.
\end{equation}
For $k=1$ this follows from the obvious identities
\begin{equation} \label{eq:x-A}
x_j A_j^{(x)} = Id + A_j^{(x)}x_j
\quad\mbox{and}\quad
x_jA_j^{(y)}= A_j^{(y)}x_j.
\end{equation}
In general it follows readily by induction.

We now proceed by induction on $r$.
Suppose (\ref{eq:comm1}) holds for some $r\ge 1$ and all $k\ge 1$.
Then
\begin{align*}
 & x_j^{r+1}  \left( A_j^{(x)}- A_j^{(y)} \right)^k = x_j
       \sum_{i=0}^r \binom{r}{i} \frac{k!}{(k-i)!}
          \left(A_j^{(x)} - A_j^{(y)}\right)^{k-i} x_j^{r-i} \\
& \; = \sum_{i=0}^r \binom{r}{i} \frac{k!}{(k-i)!}  \left[
       \left(A_j^{(x)} - A_j^{(y)}\right)^{k-i} x_j^{r+1-i}
          + (k-i) \left(A_j^{(x)} - A_j^{(y)}\right)^{k-i-1} x_j^{r-i} \right] \\
& \; =  \sum_{i=0}^{r+1} \left[ \binom{r}{i}+ \binom{r}{i-1} \right]
       \frac{k!}{(k-i)!}
            \left( A_j^{(x)}- A_j^{(y)} \right)^{k-i} x_j^{r+1-i},
\end{align*}
which completes the induction step as $\binom{r}{i}+ \binom{r}{i-1} = \binom{r+1}{i}$.
Thus (\ref{eq:comm1}) is established.

To proof (\ref{eq:comm2}), we start from
$$
(x_j - y_j)^k A_j^{(x)} = k (x_j-y_j)^{k-1} + A_j^{(x)} (x_j-y_j)^k.
$$
For $k=1$ this identity follows from (\ref{eq:x-A}) and, in general, by induction.
Finally, one proves (\ref{eq:comm2}) by induction on $s$ similarly as above.
We omit the details.
$\qed$

The next lemma is instrumental in the proof of Theorem~\ref{thm:derivative}
in the case $|\alpha|>0$.


\begin{lem} \label{lem:xAL}
If $\alpha,\beta \in \NN_0^d$ and $k\ge 1$, then
$$
  \Big| \left(A^{(x)}\right)^\alpha x^\beta  \LL_n(x,y) \Big| \le
             c_k \frac{n^{\frac{|\alpha|+|\beta|}{2}}
                   \left[K_{n+[vn]+|\alpha|+|\beta|+k}(x,x)\right]^{\frac12}
                \left[K_{n+[vn]+k}(y,y)\right]^{\frac12}}{(1+\sqrt{n} |x-y|)^k}.
$$
\end{lem}

\noindent
{\bf Proof.}
We first show that for $1 \le i \le d$
\begin{align}\label{eq:x-y2}
  & \Big|(x_i-y_i)^k \left(A^{(x)}\right)^\alpha x^\beta \LL_n(x,y)\Big| \\
  &  \qquad        \le c_k n^{(-k +|\alpha|+|\beta|)/2}
             \left[K_{n+[vn]+|\alpha|+|\beta|+k}(x,x)\right]^{\frac12}
               \left[K_{n+[vn]+k}(y,y)\right]^{\frac12}.\notag
\end{align}
Clearly
$\left(A^{(x)}\right)^{\alpha} x^\beta =
\left(A^{(x)}\right)^{\alpha-\alpha_ie_i} x^{\beta-\beta_ie_i} \cdot
  \left(A_i^{(x)}\right)^{\alpha_i} x_i^{\beta_i}
$
and the two operators separated by a dot commute.
Using (\ref{eq:comm2}) and Lemma~\ref{lem:main}, we have
\begin{align} \label{x-yAxL}
 & 2^k (x_i- y_i)^k \left(A^{(x)}\right)^{\alpha} x^{\beta}  \LL_n(x,y) \
  =
 2^k\left(A^{(x)}\right)^{\alpha-\alpha_i e_i} x^{\beta-\beta e_i} \\
& \qquad\quad \times \sum_{j=0}^{\alpha_i} \binom{\alpha_i}{j} \frac{k!}{(k-j)!}
       \left(A_i^{(x)}\right)^{\alpha_i-j}  x_i^{\beta_i}  (x_i - y_i)^{k-j} \LL_n(x,y) \notag \\
&\quad  =  \sum_{j=0}^{\alpha_i} \binom{\alpha_i}{j} \frac{2^jk!}{(k-j)!}
  \sum_{(k-j)/2 \le l \le k-j} c_{l,k-j}
      \sum_{\nu=0}^{\infty} \Delta^j \ha\Big(\frac{\nu}{n}\Big) \notag \\
&\qquad\qquad\qquad \times   \left(A^{(x)}\right)^{\alpha-je_i}  x^{\beta}
            \left(A_i^{(y)} - A_i^{(x)}\right)^{2l-k+j} \Ph_\nu(x,y). \notag
\end{align}
Furthermore, by (\ref{eq:comm1}),
\begin{align*}
\left(A_i^{(x)}\right)^{\alpha_i-j} x_i^{\beta_i} &
\left(A_i^{(y)} - A_i^{(x)}\right)^{2l-k+j}
= \sum_{\mu=0}^{\beta_i} \binom{\beta_i}{\mu} \frac{(-1)^{\mu}(2l-k+j)!}{(2l-k+j-\mu)!} \\
& \qquad\quad  \times  \left(A_i^{(x)}\right)^{\alpha_i-j}
\left(A_i^{(y)} - A_i^{(x)}\right)^{2l-k+j-\mu} x_i^{\beta_i-\mu}.
\end{align*}
As $\left(A^{(x)}\right)^{\alpha-\alpha_ie_i} x^{\beta-\beta_ie_i}$ commutes
with $A_i^{(x)}$ and $x_i$, we then conclude that
\begin{align*}
 & \left(A^{(x)}\right)^{\alpha-je_i}  x^{\beta}
            \left(A_i^{(y)} - A_i^{(x)}\right)^{2l-k+j} \Ph_\nu(x,y)
= \sum_{\mu=0}^{\beta_i} \binom{\beta_i}{\mu} \frac{(-1)^{\mu}(2l-k+j)!}{(2l-k+j-\mu)!}  \\
 & \qquad\qquad\qquad\qquad\quad \times
  \left(A^{(x)}\right)^{\alpha -je_i} \left(A_i^{(y)} - A_i^{(x)}\right)^{2l-k+j-\mu}
         x^{\beta-\mu e_i} \Ph_\nu(x,y).
\end{align*}
Using relation (\ref{eq:xPhi}) repeatedly, it follows readily that
$$
 x_i^r \Ph_\up(x) = \sum_{m=0}^r b_{m,r}(\up_i)
         \Ph_{\up+ (r-2 m)e_i}(x),
$$
where $b_{m,r}(\up_i)$ are positive numbers satisfying
$b_{m,r}(\up_i) \sim \up_i^{r/2}$.
Applying this identity to all variables we obtain
\begin{align}\label{eq:x^rPhi}
x^{\beta-\mu e_i} \Ph_\nu(x,y) =  \sum_{\mm_1=0}^{\gamma_1}\ldots
 \sum_{\mm_d=0}^{\gamma_d}  b_{\mm,\gamma}(\up)
   \sum_{|\up|=\nu} \Ph_{\up+ \gamma-2 \mm}(x)
     \Ph_{\up}(y),
\end{align}
where $\gamma_j=\beta_j$ for $j \ne i$ and
$\gamma_i = \beta_i-\mu$,
and $b_{\mm,\gamma}(\up) =b_{\mm_1,\gamma_1}(\up_1)\ldots
      b_{\mm_d,\gamma_d}(\up_{d})$.
Clearly, $|b_{\mm,\gamma}(\up)| \le c \nu^{(|\beta| -\mu)/2}$.

We now use the binomial formula and (\ref{AjHH}) to obtain
\begin{align*}
&\Big | \left(A^{(x)}\right)^{\alpha -je_i}
\left(A_i^{(y)} - A_i^{(x)}\right)^{2l-k+j-\mu}
\Ph_{\up+ \gamma-2 \mm}(x)\Ph_{\up}(y) \Big|\\
&\le c\sum_{q=0}^{2l-k+j-\mu}
\Big|\left(A_i^{(y)}\right)^{2l-k+j-\mu-q}\left(A^{(x)}\right)^{\alpha -je_i+qe_i}
\Ph_{\lambda+\gamma-2\omega}(x)\Ph_{\lambda}(y)\Big|\\
&\le c\sum_{q=0}^{2l-k+j-\mu}
\nu^{(|\alpha|+2l-k-\mu)/2}
|\Ph_{\lambda+\gamma-2\omega-je_i+qe_i}(x)\Ph_{\lambda+(2l-k+j-\mu-q)e_i}(y)\Big|
\end{align*}
and hence
\begin{align*}
&  \Big | \left(A^{(x)}\right)^{\alpha -je_i}  x^{\beta}
    \left(A_i^{(y)} - A_i^{(x)}\right)^{2l-k+j-\mu}
        \Ph_\nu(x,y) \Big|
      \le c \sum_{\mu=0}^{\beta_i}\nu^{(2l -k + |\alpha| +|\beta|-2\mu)/2} \\
&\quad  \times  \sum_{\mm_1=0}^{\gamma_1}\cdots \sum_{\mm_d=0}^{\gamma_d}
\sum_{q=0}^{2l-k+j-\mu}\sum_{|\up| = \nu}
    \left| \Ph_{\up+ \gamma - 2\mm+\alpha -je_i+q e_i}(x)
    \Ph_{\up + (2l-k+j-\mu-q)e_i}(y) \right|.
\end{align*}
As before $\Delta^l \ha(\frac{\nu}{n})=0$ if $0\le \nu\le un-l$ or $\nu\ge n+vn$,
where $0<u\le 1$ and $v>0$.
Also, by (\ref{eq:hat-a}) $|\Delta^l\ha(\frac{\nu}{n})|\le n^{-l}\|\ha^{(j)}\|_\infty$.
We use all of the above to conclude that
\begin{align*} 
& \Big| (x_i - y_i)^k \left(A_j^{(x)}\right)^\alpha
   x_j^\beta \LL_n(x,y) \Big|  \le c  n^{(-k +  |\alpha|+|\beta|)/2}
     \sum_{j=0}^{\alpha_i}  \sum_{(k-j)/2 \le l\le k-j}  |c_{l,k}|  \\
& \qquad\qquad  \times  \sum_{\nu =[un]-l}^{n+[vn]}
     \sum_{\mu=0}^{\beta_i}
         \sum_{\mm_1=0}^{\gamma_1}\cdots\sum_{\mm_d=0}^{\gamma_d}
         \left[ \sum_{q=0}^{2l-k+j-\mu}
       \Ph_{\nu+ |\alpha|+|\gamma-2\mm|+q-j}(x,x) \right]^{\frac12} \\
   &  \qquad\qquad\times
        \left[\sum_{q=0}^{2l-k+j-\mu} \Ph_{\nu+2l-k+j-\mu-q}(y,y)\right]^{\frac12}
        \notag \\
    & \qquad \le  c  n^{(-k +|\alpha|+|\beta|)/2}
         \left[K_{n+[vn]+|\alpha|+|\beta|+k}(x,x)\right]^{\frac12}
              \left[K_{n+[vn]+k}(y,y)\right]^{\frac12},    \notag
\end{align*}
where we again used the Cuachy-Schwarz inequality. This proves (\ref{eq:x-y2}).

On the other hand, using (\ref{eq:x^rPhi}) with $\mu=0$  and (\ref{AjHH}) we
can write
\begin{align*}
& \left(A^{(x)}\right)^\alpha x^\beta \LL_n(x,y) =
  \left(A^{(x)}\right)^\alpha \sum_{\nu=0}^{n+[vn]}
  \wh a \left( \frac{\nu}{n}\right)\\
&  \qquad \times \sum_{|\up|=\nu}  \sum_{\mm_1=0}^{\beta_1}\ldots
 \sum_{\mm_d=0}^{\beta_d}  b_{\mm,\beta}(\up) c_\alpha(\up)
     | \Ph_{\up+ \alpha+\beta-2 \mm}(x)  \Ph_{\up}(y)|,
\end{align*}
where $c_\alpha(\up) \sim |\up|^{|\alpha|/2}$.
Hence, using the fact that $b_{\mm,\beta}(\up) \sim |\up|^{|\beta|/2}$
we conclude that
\begin{align*}
 \Big| \left(A^{(x)}\right)^\alpha x^\beta  \LL_n(x,y) \Big|
   & \le  c \sum_{\nu=0}^{n+[vn]} \wh a  \left( \frac{\nu}{n}\right)
      \nu^{( |\alpha|+|\beta|)/2} \\
   & \qquad \times    \sum_{|\up|=\nu}
  \sum_{\mm_1=0}^{\beta_1}\ldots \sum_{\mm_d=0}^{\beta_d}
      |\Ph_{\up+ \alpha+\beta-2 \mm}(x)  \Ph_{\up}(y)| \\
   & \le c n^{( |\alpha|+|\beta|)/2} \left[K_{n+[vn]+ |\alpha|+|\beta|}(x,x)\right]^{\frac12}
              \left[K_{n+[vn]}(y,y)\right]^{\frac12}.
\end{align*}
This along with (\ref{eq:x-y2}) completes the proof of Lemma~\ref{lem:xAL}.
$\qed$

\medskip

The last step in the proof of Theorem~\ref{thm:derivative}
is to show that the operator $\partial^\alpha$ can be represented in the form
\begin{equation} \label{rep-dr}
\partial^\alpha = \sum_{\beta+\gamma\le \alpha} c_{\beta\gamma}A^\beta x^\gamma,
\end{equation}
where $\beta+\gamma\le \alpha$ means $\beta_j+\gamma_j\le \alpha_j$ for $1\le j \le d$,
and $c_{\beta\gamma}$ are constants (depending only on $\alpha$, $\beta$, $\gamma$).

By (\ref{rep-T}) $\partial_j= x_j-A_j$ and hence $\partial_j^r= (x_j-A_j)^r$.
The operators $x_j$ (multiplication by $x_j$) and $A_j$ do not commute, but
it is easy to see that
$
x_j^sA_j=sx_j^{s-1}+A_jx_j^s.
$
Applying this repeatedly one finds the representation
$$
\partial_j^r =\sum_{0\le \nu+\mu\le r}c_{\nu\mu}A_j^\nu x_j^\mu.
$$
Since the operator $A_j^\nu x_j^\mu$ commutes with $A_i^s x_i^\ell$ if $j \ne i$,
this readily implies representation (\ref{rep-dr}).

Evidently, Lemma~\ref{lem:xAL} and (\ref{rep-dr}) yield (\ref{eq:derivative1})
whenever $|\alpha|>0$.
The proof of Theorem~\ref{thm:derivative} is complete.
$\qed$

\medskip

\noindent
{\bf Proof of Lemma~\ref{lem:lower-bound}.}
Observe first that it suffices to prove (\ref{lowerbd}) only for $n$ sufficiently
large since it holds trivially if $2/\del \le n \le c$.

We next prove (\ref{lowerbd}) for $d =1$. The Christoph-Darboux formula for
the Hermite polynomials (\cite[(5.59)]{Sz})
shows that
$$
 K_m(x,x) = (2^{m+1} m!)^{-1}
     \left[H'_{m+1}(x) H_m(x) - H_m'(x) H_{m+1}(x)\right].
$$
Using the fact that $H_{m+1}'(x) = 2(m+1) H_m(x)$ (\cite[(5.5.10)]{Sz})
and $H_{m+1}(x) = 2 x H_m(x) - 2m H_{m-1}(x)$ (\cite[(5.5.8)]{Sz}), we can
rewrite $K_m(x,x)$ as
$$
  K_m(x,x) = (2^{m+1} m!)^{-1} \left[2 (m+1) H_m^2 (x)
         - 4 m x H_{m-1}(x) H_m(x) + 2 m^2 H_{m-1}^2 (x)\right].
$$
Written in terms of the orthonormal Hermite functions, the above identity
becomes
$$
  (m+1) h_m^2(x) +  m h_{m-1}^2(x) = K_m(x,x) + \sqrt{2m}\,
     x \,h_m(x) h_{m-1}(x).
$$
In particular, it follows that for $|x| \le 2 \sqrt{2m+1}$
$$
  (m+1) h_m^2(x) +  m h_{m-1}^2(x) \ge  K_m(x,x) - 2 \sqrt{2m}
         \sqrt{2m+1} |h_m(x)| \cdot |h_{m-1}(x)|.
$$
and hence
\begin{align*}
 (3 m + 2) h_m^2(x) +  3 m h_{m-1}^2(x) & \ge  K_m(x,x) +
           \left (\sqrt{2m+1} |h_m(x)| - \sqrt{2m} |h_{m-1}(x)|\right)^2 \\
     & \ge K_m(x,x).
\end{align*}
Consequently, for $|x| \le 2 \sqrt{2n+1}$
\begin{align*}
\sum_{m=[(1-\del)n]}^n h_m^2(x)
      & \ge \frac{1}{2} \sum_{m=[(1-\del)n]}^{n}
         \left(h_m^2(x) + h_{m+1}^2(x)\right) \\
    & \ge \frac{c}{n}  \sum_{m=[(1-\del)n]}^n K_m(x,x)
     \ge c_1\, K_{[\rho n]}(x,x),
\end{align*}
which proves (\ref{lowerbd}) when $d =1$.

\smallskip

For $d > 1$ we need the following identity which follows from the generating
function of Hermite polynomials (see e.g. \cite{Th}):
$$
  \sum_{k=0}^\infty \CH_k(x,x) r^k
     = \pi^{-d/2} (1-r^2)^{-d/2} e^{- \frac{1-r}{1+r} \|x\|^2} : = F_d(r,t),
$$
where $t = |x|$.
Let us denote $\CH_{k,d}(x,x)=\CH_k(x,x)$ for $x \in \RR^d$
in order to indicate the dependence on $d$.
Then it follows from above that
\begin{align} \label{recursive}
\sum_{k=0}^\infty r^k [\CH_{k,d}(x,x) - \CH_{k-2,d}(x,x)] & =
  (1-r^2) \sum_{k=0}^\infty r^k \CH_{k,d}(x,x)\\
      & = (1-r^2) F_d(r,t) = \pi^{-1} F_{d-2}(r,t). \notag
\end{align}
Notice that $\CH_{n,d}(x,x)$ is a radial function and hence a function of $t$.
Thus comparing the coefficients of $r^k$ in both side shows that
$$
  \CH_{k,d}(x,x) -   \CH_{k-2,d}(x,x) = \pi^{-1} \CH_{k,d-2}(x,x),
$$
which implies
$$
 \CH_{k,d}(x,x) + \CH_{k-1,d}(x,x)= \pi^{-1} \sum_{j=0}^k \CH_{j,d-2}(x,x)
      = \pi^{-1}K_{k,d-2}(x,x).
$$
Now, summing over $k$ we get
\begin{align*}
 & \sum_{k= [(1-\del)n]}^{n} \CH_{k,d}(x,x) \ge \frac{1}{2}
    \sum_{k=[(1-\del)n]+1}^{n}
      [ \CH_{k,d}(x,x)+\CH_{k-1,d}(x,x) ] \\
 &     \ge c   \sum_{k=[(1-\del)n]+1}^{n} K_{k,d-2} (x,x)
       \ge c\,n  K_{[(1-\eps)n], d-2} (x,x) \ge c \, n
         \sum_{k= [(1-\del)n]}^{[(1-\eps)n]} \CH_{k,d-2}(x,x),
\end{align*}
where $\eps:=(1-\rho)/d$.
Evidently, by induction this estimate yields (\ref{lowerbd}) for $d$ odd.

\smallskip

To establish the result for $d$ even, we only have to prove estimate (\ref{lowerbd})
for $d = 2$. By the definion of $F_d(r,t)$, we have
\begin{align*}
 F_0(r,t) & = e^{-\frac{1-r}{1+r} t^2} = \pi^{1/2}(1-r^2)^{1/2} F_1(r,t) \\
   & = \pi^{1/2} \sum_{j=0}^\infty \binom{1/2}{j} (-1)^j r^{2j}
          \sum_{n=0}^\infty h_k^2(t) r^k
   = \pi^{1/2} \sum_{k=0}^\infty \left[\sum_{j=0}^k
              a_{k-j} h_j^2(x) \right] r^k,
\end{align*}
where $a_{2j} = (-1)^j \binom{1/2}{j}$ and $a_{2j-1} =0$.
Hence, using (\ref{recursive})
$$
\CH_{k,2}(x,x) - \CH_{k-2,2}(x,x) = \pi^{-1/2} \sum_{j=0}^k a_{k-j} h_j^2(t).
$$
Consequently,
\begin{align*}
 \CH_{k,2}(x,x) + & \CH_{k-1,2}(x,x)  = \pi^{-1/2}
  \sum_{l=0}^k \sum_{j=0}^l a_{l-j} h_j^2(t)
    = \sum_{j=0}^k h_j^2(t) \sum_{l=0}^{k-j} a_l.
\end{align*}
A simple combinatorial formula shows that
$$
\sum_{l=0}^{k-j} a_l  = \sum_{l=0}^{[(k-j)/2]} (-1)^l \binom{1/2}{l}
    =  \frac{\Gamma(\frac{1}{2} + [\frac{k-j}{2}])}
        {\Gamma(\frac{1}{2})\Gamma(1+ [\frac{k-j}{2}])},
$$
which is positive for all $0 \le j \le k$.
Furthermore, by
$\Gamma(k+a)/\Gamma(k+1) \sim k^{a-1}$ it follows that
$ \sum_{l=0}^{k-j} a_l \ge c k^{-1/2}$ for $0 \le j \le \alpha k$
for any $\alpha <1$.
Therefore,
$$
  \CH_{k,2}(x,x) +  \CH_{k-1,2}(x,x)  \ge c k^{-1/2}
     \sum_{j=0}^{\alpha k} h_j^2(t)
$$
and summing over $k$ we get
$$
  \sum_{k= [(1-\del)n]}^{n} \CH_{k,2}(x,x) \ge c n^{1/2}
       K_{[\rho n],1}(t,t),
$$
which establishes (\ref{lowerbd}) for $d =2$.
$\qed$

\subsection{Proofs for Sections~\ref{Tri-Liz-spaces}-\ref{Besov-spaces}}
\label{proofs1}
${}$

\medskip

\noindent
{\bf Proof of Lemma \ref{lem:sum<M}.}
Let
\begin{equation}\label{def-b-diam}
b_j^\dm(x):= \sum_{\eta\in\cX_j}\frac{|b_\eta|}{(1+2^jd(x, R_\eta))^\sigma},
\end{equation}
where $d(x, E)$ stands for the $\ell^\infty$ distance of $x$ from $E\subset \R^d$.
Evidently,
\begin{equation}\label{bj*<bj-dm}
b_j^*(x)\le cb_j^\dm(x) \quad\mbox{and}\quad
b_\xi^*\ONE_{R_\xi}(x) \le cb_j^\dm(x), \quad x\in \R^d, \quad \xi\in\cX_j.
\end{equation}
We will show that
\begin{equation}\label{bj-dm<M}
b_j^\dm(x) \le c \cM_s\Big(\sum_{\omega\in\cX_j}|b_\omega|\ONE_{R_\omega}\Big)(x), \quad x\in \R^d.
\end{equation}
In view of (\ref{bj*<bj-dm}) this implies (\ref{s*<M})-(\ref{sum<M}),
and hence Lemma~\ref{lem:sum<M}.

By the construction of the tiles $\{R_\xi\}$ in (\ref{def.Q-xi})-(\ref{def.Q-j})
it follows that there exists a constant $\cd>0$ depending only on $d$
such that
$$
Q_j:=\cup_{\xi\in\cX_j} R_\xi \subset [-\cd 2^j, \cd 2^j]^d.
$$

Fix $x\in \R^d$. To prove (\ref{bj-dm<M}) we consider two cases for $x$.

\medskip

{\em Case 1.} $|x|_\infty>2\cd 2^j$.
Then $d(x, R_\eta)>|x|_\infty/2$ for $\eta\in\cX_j$ and hence
\begin{align}\label{est-b*}
b_j^\dm(x)
= \sum_{\eta\in\cX_j}\frac{|b_\eta|}{(1+2^jd(x, R_\eta))^\sigma}
&\le \frac{c}{(2^{j}|x|_\infty)^\sigma}\sum_{\eta\in\cX_j}|b_\eta| \notag\\
&\le \frac{c4^{jd\lambda}}{(2^{j}|x|_\infty)^\sigma}
\Big(\sum_{\eta\in\cX_j}|b_\eta|^s\Big)^{1/s},
\end{align}
where $\lambda:=1-\min\{1, 1/s\}$ and for the last estimate
we use H\"{o}lder's inequality if $s>1$ and the $s$-triangle inequality
if $s < 1$.

Denote $Q_x:=[-|x|_\infty, |x|_\infty]^d$. Notice that $Q_j\subset Q_x$.
From above we infer
\begin{align*}
b_j^\dm(x)
&\le \frac{c4^{jd\lambda}|x|_\infty^{d/s}}{(2^{j}|x|_\infty)^\sigma}
\Big(\frac{1}{|Q_x|}\int_{Q_x}
\Big(\sum_{\eta\in\cX_j}|b_\eta|\ONE_{R_\xi}(y)\Big)^s dy\Big)^{1/s}\\
&\le c2^{j(2d\lambda-\sigma)}|x|^{d/s-\sigma}
\cM_s\Big(\sum_{\eta\in\cX_j}|b_\eta|\ONE_{R_\xi}\Big)(x)
\le c\cM_s\Big(\sum_{\eta\in\cX_j}|b_\eta|\ONE_{R_\xi}\Big)(x)
\end{align*}
as claimed.
Here we used the fact that $\sigma\ge d\max\{2, 1/s\}$.

\medskip

{\em Case 2.} $|x|_\infty\le 2\cd 2^j$.
To make the argument more transparent we first subdivide the tiles
$\{R_\eta\}_{\eta\in\cX_j}$ into boxes of almost equal sides of length
$\sim 2^{-j}$.
By the construction of the tiles (see (\ref{def.Q-xi})) there exists
a constant $\tilde{c}>0$ such that the minimum side of each tile $R_\eta$
is $\ge \tilde{c}2^{-j}$.
Now, evidently each tile $R_\eta$ can be subdivided into a disjoint
union of boxes $R_\theta$ with centers $\theta$ such that
$$
\theta+[-\tilde{c}2^{-j-1}, \tilde{c}2^{-j-1}]
\subset R_\theta \subset
\theta+[-\tilde{c}2^{-j}, \tilde{c}2^{-j}].
$$
Denote by $\widehat{\cX}_j$ the set of centers of all boxes obtained
by subdividing the tiles from $\cX_j$.
Also, set $b_\theta:=b_\eta$ if $R_\theta\subset R_\eta$.
Evidently,
\begin{equation}\label{est-b*2}
b_j^\dm(x)
:= \sum_{\eta\in\cX_j}\frac{|b_\eta|}{(1+2^jd(x, R_\eta))^\sigma}
\le \sum_{\theta\in\hcX_j}\frac{|b_\theta|}{(1+2^jd(x, R_\theta))^\sigma}
\end{equation}
and
\begin{equation}\label{est-b-b}
\sum_{\eta\in \cX_j}|b_\eta|\ONE_{R_\eta}
= \sum_{\eta\in \hcX_j}|b_\theta|\ONE_{R_\theta}.
\end{equation}

Denote
$Y_0:=\{\theta\in\hcX_j: 2^j|\theta-x|_\infty\le \tilde{c}\}$,
\begin{align*}
&Y_m:=\{\theta\in\hcX_j: \tilde{c}2^{m-1}\le 2^j|\theta-x|_\infty\le \tilde{c}2^m\},
\quad \mbox{and}\\
&Q_m:=\{y\in\R^d: |y-x|_\infty\le \tilde{c}(2^m+1)2^{-j}\}, \quad m\ge 1.
\end{align*}
%
%
Clearly,
$\#Y_m \le c2^{md}$,
$\cup_{\theta\in Y_m}R_\theta \subset Q_m$, and
$\hcX=\cup_{m\ge 0} Y_m$.
Similarly as in (\ref{est-b*})
\begin{align*}
\sum_{\theta\in Y_m}\frac{|b_\theta|}{(1+2^jd(x, R_\theta))^\sigma}
&\le c2^{-m\sigma}\sum_{\theta\in Y_m}|b_\theta|
\le c2^{-m\sigma}2^{md\lambda}
\Big(\sum_{\theta\in Y_m}|b_\theta|^s\Big)^{1/s}\\
&\le c2^{-m(\sigma-d\lambda-d/s)}
\Big(\frac{1}{|Q_m|}\int_{Q_m}
\Big(\sum_{\theta\in Y_m}|b_\theta|\ONE_{R_\theta}(y)\Big)^s dy\Big)^{1/s}\\
&\le c2^{-m(\sigma-d\max\{1, 1/s\})}
\cM_s\Big(\sum_{\eta\in \cX_j}|b_\eta|\ONE_{R_\eta}\Big)(x),
\end{align*}
where we used (\ref{est-b-b}).
Summing up over $m\ge 0$, taking into account that
$\sigma > d\max\{2, 1/s\})$, and also using (\ref{est-b*2})
we arrive at (\ref{bj-dm<M}).
$\qed$


\medskip

\noindent
{\bf Proof of Lemma \ref{lem:a*=b*}.}
For this proof we will need an additional lemma.


\begin{lem}\label{lem:g-g}
Suppose $g \in \V_{4^j}$ and $\xi \in\cX_j$.
Then for any $k>0$ and $L>0$
we have for $x', x''\in 2R_\xi$
\begin{equation}\label{g-g1}
|g(x')-g(x'')|\le c2^j|x'-x''|\sum_{\eta \in \cX_j}
\frac{|g(\eta)|}{(1+2^j|\xi-\eta|)^k}
\end{equation}
and
\begin{equation}\label{g-g2}
|g(x')-g(x'')|\le \hat{c}2^{-jL}|x'-x''|\sum_{\eta \in \cX_j}
\frac{|g(\eta)|}{(1+2^j|\xi-\eta|)^k},
\; \mbox{if}\; |\xi|_\infty > (1+2\delta)2^{j+1}.
\end{equation}
Here
$c$, $\hat{c}$ depend on $k$, $d$, and $\delta$,
and $\hat{c}$ depends on $L$ as well;
$2R_\xi\subset \RR^d$ is the set obtained by dilating $R_\xi$
by a factor of 2 and with the same center.
\end{lem}

\noindent
{\bf Proof.}
Let $\Lambda_{4^j}$ be the kernel from (\ref{def-LLn})
with $n=4^j$, where $\ha$ is admissible of type (a)
with $v:= \delta$.
Then
$\Lambda_{4^j}*g=g$, since $g \in \V_{4^j}$,
and  $\Lambda_{4^j}(x, \cdot)\in V_{[(1+\delta)4^j]}$.
Note that $[(1+\delta)4^j]+4^j \le 2N_j-1$.
Therefore,
we can use the cubature formula from Corollary~\ref{cor:cubature}
to obtain
$$
g(x)=\int_{\R^d} \Lambda_{4^j}(x, y)g(y)dy
= \sum_{\eta\in\cX_j}\lambda_\eta\Lambda_{4^j}(x, \eta)g(\eta),
$$
where the weights $\lambda_\xi$ obey (\ref{tiles})
(see also (\ref{tiles1})-(\ref{tiles2})).
Hence, for $x', x''\in 2R_\xi$
\begin{align}\label{est-g-g}
|g(x')-g(x'')|
&\le \sum_{\eta\in\cX_j}\lambda_\eta
|\Lambda_{4^j}(x', \eta)-\Lambda_{4^j}(x'', \eta)||g(\eta)|\notag\\
&\le c|x'-x''|\sum_{\eta\in\cX_j}\lambda_\eta\sup_{x\in 2R_\xi}
|\nabla \Lambda_{4^j}(x, \eta)||g(\eta)|.
\end{align}

Note that
$\Big(4([(1+\delta)4^j]+k+1)+2\Big)^{1/2} \le (1+\delta)2^{j+1}$
for sufficiently large $j$ (depending on $k$ and $\delta$).
Therefore, we have from (\ref{est:derivative2})-(\ref{est:derivative3})
\begin{equation}\label{est-nabla1}
|\nabla \Lambda_{4^j}(x, \eta)|
\le \frac{c2^{j(d+1)}}{(1+2^j|x-\eta|)^k},
\quad x\in\RR^d, \; \eta\in\cX_j,
\end{equation}
and for any $L>0$ (we need $L\ge k$)
\begin{equation}\label{est-nabla2}
|\nabla \Lambda_{4^j}(x, \eta)|
\le \frac{c2^{-2jL}}{(1+2^j|x-\eta|)^k},
\;\mbox{if} \; |x|_\infty > (1+\delta)2^{j+1}\;
\mbox{or} \; |\eta|_\infty > (1+\delta)2^{j+1}.
\end{equation}

Suppose $|\xi|_\infty > (1+2\delta)2^{j+1}$, then
$
2R_\xi \subset \{x\in\RR^d: |x|_\infty > (1+\delta)2^{j+1}\}
$
for sufficiently large $j$.
Combining (\ref{est-g-g}) with (\ref{est-nabla2})
and (\ref{tiles2}) we get
\begin{equation}\label{g-g10}
|g(x')-g(x'')|
\le c2^{-jd/3}|x'-x''|\sum_{\eta\in\cX_j}
\sup_{x\in 2R_\xi}\frac{2^{-j(k+L)}}{(1+2^j|x-\eta|)^k}|g(\eta)|,
\end{equation}
where we used that ${\rm diam}\, (2R_\xi) \le c 2^{-j/3}$.
However, for any $x\in 2R_\xi$ we have
$$
1+2^j|\xi-\eta| \le 1+2^j(|\xi-x| +|x-\eta|) \le 1+2^j(c2^{-j/3} +|x-\eta|)
\le c2^j(1+2^j|x-\eta|).
$$
We use this in (\ref{g-g10}) to obtain (\ref{g-g2})
for sufficiently large $j$.

One proves (\ref{g-g1}) in a similar fashion.
In the case $j\le c$ estimates (\ref{g-g1})-(\ref{g-g2})
follow easily by (\ref{est-g-g}).
$\qed$


\medskip

We are now prepared to prove Lemma~\ref{lem:a*=b*}.
Let $g\in V_{4^j}$.
Pick $\ell \ge 1$ sufficiently large (to be determined later on) and denote
for $\xi\in\cX_j$
\begin{equation}\label{def-Xj-xi}
\cX_{j+\ell}(\xi):=\{\eta\in\cX_{j+\ell}: R_\eta\cap R_\xi \ne \emptyset\}
\quad \mbox{and}
\end{equation}
\begin{equation}\label{def-d-xi}
d_\xi:=\sup\{|g(x')-g(x'')|:
x',x''\in R_\eta \;\; \mbox{for some} \;\;\eta\in\cX_{j+\ell}(\xi)\}.
\end{equation}

We first estimate $d_\xi$, $\xi\in\cX_j$.


{\em Case A:} $|\xi|_\infty\le (1+3\delta)2^{j+1}$.
By (\ref{tiles1}) it follows that for sufficiently large $\ell$ (depending only on $d$ and $\delta$)
$\cup_{\eta\in\cX_{j+\ell}(\xi)} R_\eta \subset 2R_\xi$.
Hence, using Lemma~\ref{lem:g-g} (see (\ref{g-g1})) with $k\ge \sigma$, we get
\begin{equation}\label{est-d-xi1}
d_\xi \le c2^{-\ell}\sum_{\eta\in\cX_j}\frac{|g(\eta)|}{(1+2^j|\xi-\eta|)^{\sigma}},
\end{equation}
for sufficiently large $j$ (depending only on $d$ and $\delta$),
where $c>0$ is a constant independent of $\ell$.


{\em Case B:} $|\xi|_\infty>(1+3\delta)2^{j+1}$.
By (\ref{tiles1})
$|x|_\infty > (1+2\delta)2^{j+1}$ for
$x\in \cup_{\eta\in\cX_{j+\ell}(\xi)} R_\eta$
if $j$ is sufficiently large.
We apply estimate (\ref{g-g2}) of Lemma~\ref{lem:g-g} with $k\ge \sigma$ and
$L=1$ to obtain
\begin{equation}\label{est-d-xi2}
d_\xi \le c2^{-j}\sum_{\eta\in\cX_j}\frac{|g(\eta)|}{(1+2^j|\xi-\eta|)^{\sigma}}.
\end{equation}

To estimate $M_\xi^*$, $\xi\in\cX_j$, we consider two cases for $\xi$.

\smallskip

{\em Case 1:} $|\xi|_\infty\le (1+4\delta)2^{j+1}$.
By (\ref{tiles1}), we have for sufficiently large $j$:
\begin{equation}\label{incl-R}
R_\xi\sim \xi+[-2^{-j}, 2^{-j}]^d,
\quad \mbox{and}\quad
R_\eta\sim \eta+[-2^{-j-\ell}, 2^{-j-\ell}]^d,
\quad \eta\in \cX_{j+\ell}(\xi).
\end{equation}
By the definition of $d_\xi$ in (\ref{def-d-xi}) it follows that
$M_\xi \le m_\lambda+d_\xi$ for some $\lambda\in\cX_{j+\ell}(\xi)$
and hence, using (\ref{incl-R}),
$$
M_\xi \le c\sum_{\omega\in\cX_{j+\ell}}\frac{m_\omega}{(1+2^{j+\ell}|\xi-\omega|)^\sigma}
+ d_\xi =: G_\xi+ d_\xi,
\quad c=c(d, \delta, \sigma, \ell).
$$
Consequently,
\begin{equation}\label{Est-M*}
M_\xi^* \le G_\xi^*+ d_\xi^*.
\end{equation}
Write $\cX_j':= \{\eta\in\cX_j: |\eta|_\infty\le (1+3\delta)2^{j+1}\}$
and
$\cX_j'':= \cX_j\setminus \cX_j'$.
Now, we use (\ref{est-d-xi1})-(\ref{est-d-xi2}) to obtain
\begin{align*}
d_\xi^*
=\sum_{\eta\in\cX_j}\frac{d_\eta}{(1+2^j|\xi-\eta|)^\sigma}
&\le c2^{-\ell}\sum_{\eta\in\cX_j}\sum_{\omega\in\cX_j'}
\frac{|g(\omega)|}{(1+2^j|\xi-\eta|)^\sigma(1+2^j|\eta-\omega|)^\sigma}\\
&+c2^{-j}\sum_{\eta\in\cX_j}\sum_{\omega\in\cX_j''}
\frac{|g(\omega)|}{(1+2^j |\xi-\eta|)^\sigma(1+2^j|\eta-\omega|)^\sigma},
\end{align*}
replacing $\cX_j'$ and $\cX_j''$ by $\cX_j$ above and shifting the order of summation
we get
\begin{align}\label{Est-d-xi}
d_\xi^*
&\le c(2^{-\ell}+2^{-j})\sum_{\omega\in\cX_j}|g(\omega)| \sum_{\eta\in\cX_j}
\frac{1}{(1+2^j|\xi-\eta|)^\sigma(1+2^j|\eta-\omega|)^\sigma}\\
&\le c(2^{-\ell}+2^{-j})\sum_{\omega\in\cX_j}
\frac{|g(\omega)| }{(1+2^j|\xi-\omega|)^\sigma}
\le c(2^{-\ell}+2^{-j})M_\xi^*.\notag
\end{align}
Here the constant $c$ is independent of $\ell$ and $j$,
and we used that
\begin{align}
\sum_{\eta\in\cX_j}\frac{1}{(1+2^j|\xi-\eta|)^\sigma(1+2^j|\eta-\omega|)^\sigma}
&\le \int_{\R^d}\frac{c2^{jd}}{(1+2^j|\xi-y|)^\sigma(1+2^j|y-\omega|)^\sigma}dy\notag\\
&\le \frac{c}{(1+2^j|\xi-\omega|)^\sigma}
\quad (\sigma>d). \label{sum<int}
\end{align}
These estimates are standard and easy to prove utilizing the fact that the tiles
$\{R_\eta\}_{\eta\in\cX_j}$
do not overlap and obey (\ref{tiles1}).

To estimate $G_\xi^*$ we use again (\ref{tiles1}) and (\ref{sum<int}). We get
\begin{align*}
G_\xi^*
&=\sum_{\eta\in\cX_j}\frac{G_\eta}{(1+2^j|\xi-\eta|)^\sigma}
\le c\sum_{\eta\in\cX_j}\sum_{\omega\in\cX_{j+\ell}}
\frac{m_\omega}{(1+2^j|\xi-\eta|)^\sigma(1+2^j|\eta-\omega|)^\sigma}\\
&\le c\sum_{\omega\in\cX_{j+\ell}}m_\omega\sum_{\eta\in\cX_j}
\frac{1}{(1+2^j |\xi-\eta|)^\sigma(1+2^j|\eta-\omega|)^\sigma}
\le c \sum_{\omega\in\cX_{j+\ell}}
\frac{m_\omega}{(1+2^j |\xi-\omega|)^\sigma}\\
&\le c2^{\ell\sigma} \sum_{\omega\in\cX_{j+\ell}}
\frac{m_\omega}{(1+2^{j+\ell}|\lambda-\omega|)^\sigma}
=cm_\lambda^*
\quad\mbox{for each $\lambda \in \cX_{j+\ell}(\xi)$.}
\end{align*}
Combining this with (\ref{Est-M*})-(\ref{Est-d-xi}) we obtain
$$
M_\xi^*\le c_1 m_\lambda^*+ c_2(2^{-\ell}+2^{-j})M_\xi^*
\quad\mbox{for $\lambda \in \cX_{j+\ell}(\xi)$,}
$$
where $c_2>0$ is independent of $\ell$ and $j$.
Choosing $\ell$ and $j$ sufficiently large (depending only on $d$, $\delta$, and $\sigma$)
this yields $M_\xi^*\le cm_\lambda^*$ for all $\lambda \in \cX_{j+\ell}(\xi)$.
For $j\le c$ this relation follows as above but using only (\ref{g-g1})
and taking $\ell$ large enough. We skip the details.
Thus (\ref{a*=b*}) is established in Case 1.

\medskip


{\em Case 2:}  $|\xi|_\infty > (1+4\delta)2^{j+1}$.
In this case for sufficiently large $j$ (depending only on $d$, $\delta$, and $\sigma$)
$
|x|_\infty \ge (1+3\delta)2^{j+1}
$
for $x\in\cup_{\eta\in\cX_{j+\ell}(\xi)} R_\eta$.
Hence, using (\ref{g-g2}) with $L=1$, we have
$$
M_\xi\le m_\omega+c2^{-j}\sum_{\eta\in\cX_j}\frac{|g(\eta)|}{(1+2^j|\xi-\eta|)^\sigma}
\le m_\omega +c2^{-j}M_\xi^*
\quad\mbox{for all $\omega\in\cX_{j+\ell}(\xi)$,}
$$
where $c>0$ is independent of $j$.
Fix $\lambda\in\cX_{j+\ell}(\xi)$ and for each $\eta\in\cX_j$, $\eta\ne \xi$,
choose $\omega_\eta\in \cX_{j+\ell}(\eta)$ so that
$
|\lambda-\omega_\eta|=\min_{\omega\in\cX_{j+\ell}(\eta)} |\lambda-\omega|.
$
Then from above
\begin{equation}\label{Est-M*2}
M_\xi^* \le \sum_{\eta\in\cX_j}\frac{m_{\omega_\eta}}{(1+2^j|\xi-\eta|)^\sigma}
+c2^{-j}\sum_{\eta\in\cX_j}\frac{M_\eta^*}{(1+2^j|\xi-\eta|)^\sigma}
=: A_1+A_2.
\end{equation}
By (\ref{almost-eq}) it easily follows that $\omega_\eta$ from above obeys
$|\lambda-\omega_\eta| \le c|\xi-\eta|$
and hence
\begin{equation}\label{Est-A1}
A_1\le c\sum_{\eta\in\cX_j}\frac{m_{\omega_\eta}}{(1+2^j|\lambda-\omega_\eta|)^\sigma}
\le c2^{\ell\sigma}\sum_{\omega\in\cX_{j+\ell}}\frac{m_{\omega}}{(1+2^{j+\ell}|\lambda-\omega|)^\sigma}
\le c_1m_\lambda^*.
\end{equation}
On the other hand, using Definition~\ref{def-s*} and (\ref{sum<int}), we have
\begin{align*}
A_2
&\le c2^{-j}\sum_{\eta\in\cX_j}\sum_{\omega\in\cX_j}
\frac{M_\omega}{(1+2^j|\xi-\eta|)^\sigma(1+2^j|\eta-\omega|)^\sigma}\\
&\le c2^{-j}\sum_{\omega\in\cX_j} M_\omega\sum_{\eta\in\cX_j}
\frac{1}{(1+2^j|\xi-\eta|)^\sigma(1+2^j|\eta-\omega|)^\sigma}\\
&\le c_22^{-j}\sum_{\omega\in\cX_j}
\frac{M_\omega}{(1+2^j|\eta-\omega|)^\sigma}
= c_22^{-j}M_\omega^*,
\end{align*}
where $c_2>0$ is independent of $j$.
Combining this with (\ref{Est-M*2})-(\ref{Est-A1}) we arrive at
$$
M_\xi^*\le c_1m_\lambda^*+ c_22^{-j}M_\xi^*
\quad\mbox{for $\lambda\in\cX_{j+\ell}(\xi)$.}
$$
Choosing $j$ sufficiently large we get
$
M_\xi^*\le c_1m_\lambda^*
$
for each $\lambda\in\cX_{j+\ell}(\xi)$.
For $j\le c$ this estimate follows as in Case~1 but using only (\ref{g-g1}).
This completes the proof of Lemma~\ref{lem:a*=b*}.

\medskip


\noindent
{\bf Proof of Lemma~\ref{lem:FH-diff-F}.}
To prove (\ref{FH-diff-F1}) we first show that there exit
$\eps>0$ and $r>1$ such that
\begin{equation}\label{F-embed3}
F^{\a q}_p(H) \hookrightarrow F^{\eps 2}_r(H).
\end{equation}
Indeed, if $p>1$, using that $\a >0$, Proposition~\ref{prop:F-embed-F} (a) yields
$
F^{\a q}_p \hookrightarrow F^{\eps 2}_p
$
for any $0<\eps<\alpha$.
On the other hand, if $p\le 1$, then $\alpha-d/p>-d$ and hence there exist $\delta>0$ and $\eps>0$
such that, first, $\alpha-d/p>-d/(1+\delta)$ and then
$\alpha-d/p=\eps-d/(1+\delta)$.
Now, by Proposition~\ref{prop:F-embed-F} (b) we have
$
F^{\a q}_p \hookrightarrow F^{\eps 2}_{1+\delta}.
$
Thus (\ref{F-embed3}) is established.

Denote $h_y(x):=h(x-y)$.
It follows by Proposition~\ref{prop:identification} and
Theorem~\ref{thm:needlet-Tr-Liz} that
$$
\|h_y\|_r
\sim \Big\|\Big(\sum_{\xi\in\cX}
\Big(|R_\xi|^{-1/2}|\langle h_y, \ph_\xi\rangle|
\ONE_{R_\xi}(\cdot)\Big)^2\Big)^{1/2}\Big\|_r
=:\CN(h_y).
$$
Fix $J\ge 1$ and denote $\CY_J:=\cup_{0\le j\le J}\cX_j$.
By the decay of needlets (see (\ref{local-Needlets21}))
it follows that
$$
\max_{\xi\in\CY_J}|\langle h_y, \ph_\xi\rangle| \to 0
\quad\mbox{as $|y|\to \infty$.}
$$
Hence there exists $A>0$ such that if $|y|>A$,
\begin{equation}\label{YJ-sum}
\Big\|\Big( \sum_{\xi \in \CY}
\Big(|R_\xi|^{-1/2}|\langle h_y, \ph_\xi\rangle|\ONE_{R_\xi}(\cdot)\Big)^2\Big)^{1/2}\Big\|_r
\le \frac{1}{2}\CN(h_y).
\end{equation}
Evidently, $h_y$ being $C^\infty$ and compactly supported belongs to all
Hermite-F-spaces and by (\ref{F-embed3})
$\|h_y\|_{\Fapq(H)} \ge c\|h_y\|_{F^{\eps 2}_r(H)}$.
We now use Theorem~\ref{thm:needlet-Tr-Liz}
and (\ref{YJ-sum}) to obtain, for $|y|>A$,
\begin{align*}
\|h_y\|_{\Fapq(H)} &\ge c\|h_y\|_{F^{\eps 2}_r(H)}
\ge c\Big\|\Big(  \sum_{j=0}^\infty  2^{\eps j}\sum_{\xi\in\cX_j}
\Big(|R_\xi|^{-1/2}|\langle h_y, \ph_\xi\rangle|
\ONE_{R_\xi}(\cdot)\Big)^2\Big)^{1/2}\Big\|_r\\
&\ge c2^{J\eps}\Big\|\Big(\sum_{j=J+1}^\infty \sum_{\xi\in\cX\setminus\CY_J}
\Big(|R_\xi|^{-1/2}|\langle h_y, \ph_\xi\rangle|
\ONE_{R_\xi}(\cdot)\Big)^2\Big)^{1/2}\Big\|_r\\
&\ge (1/2)c 2^{J\eps}
\Big\|\Big(\sum_{\xi\in\cX}
\Big(|R_\xi|^{-1/2}|\langle h_y, \ph_\xi\rangle|
\ONE_{R_\xi}(\cdot)\Big)^2\Big)^{1/2}\Big\|_r\\
&\ge c'2^{J\eps}\|h_y\|_r
= c'2^{J \eps}\|h\|_r
\quad (\|h\|_r>0),
\end{align*}
where $c'>0$ is independent of $J$.
Letting $J\to\infty$ the above implies (\ref{FH-diff-F1}).

\smallskip

We next prove (\ref{FH-diff-F2}).
Choose $k>\max\{\alpha+d, d/p\}$.
Using (\ref{local-Needlets21})-(\ref{local-Needlets22}) we get,
for $\xi\in\cX_j$ and $|\xi-y|>|y|/2$,
and sufficiently large $|y|$,
$$
|\langle h_y, \ph_\xi\rangle| \le \frac{c2^{jd/2}}{(1+2^j|y-\xi|)^k}
\le \frac{c'2^{jd/2}}{(1+2^j|y-x|)^k}
\quad\mbox{for each $x\in R_\xi$.}
$$
Hence, using also (\ref{tiles3}) we have that for $|x-y|\ge |y|/4$ and
$|y|$ sufficiently large
\begin{align*}
G(x) &:=\sum_{j=0}^\infty  2^{j\a q}\sum_{\xi\in\cX_j, |\xi-y|>|y|/2}
\Big(|R_\xi|^{-1/2}|\langle h_y, \ph_\xi\rangle|
\ONE_{R_\xi}(x)\Big)^q\\
&\le c\sum_{j=0}^\infty
\frac{2^{j(\a+d)q}}{(1+2^j|y-x|)^{kq}}
\le \frac{c}{|y-x|^{kq}}
\sum_{j=0}^\infty
2^{-j(k-\a-d)q}
\le \frac{c}{|y-x|^{kq}},
\end{align*}
while
$$
G(x) = 0\quad \mbox{if $|x-y|<|y|/4$.}
$$
Hence,
$$
\|h_y\|_{F^*_y} \le c\Big(\int_{|x-y|>|y|/4}\frac{dx}{|y-x|^{kp}}\Big)^{1/p}
\le \frac{c}{|y|^{k-d/p}},
$$
which yields (\ref{FH-diff-F2}).
$\qed$

\medskip


\noindent
{\bf Proof of Lemma~\ref{lem:Polyn-Lp}.}
Let $g\in \V_{4^j}$ ($j\ge 0$) and $0<p<\infty$.
We will utilize Definition~\ref{def-s*} and
Lemmas~\ref{lem:sum<M}-\ref{lem:a*=b*}.
To this end choose $0<s<\min\{p, 1\}$ and $\sigma > d\max\{2, 1/s\}$.
Set $M_\xi:=\sup_{x\in R_\xi}|g(x)|$, $\xi\in\cX_j$,
and $m_\lambda:=\inf_{x\in R_\lambda}|g(x)|$, $\lambda\in\cX_{j+\ell}$,
where $\ell\ge 1$ is the constant from Lemma \ref{lem:a*=b*}.
Using Lemmas~\ref{lem:sum<M}-\ref{lem:a*=b*} and the maximal inequality (\ref{max-ineq})
we get
\begin{align*}
&\Big(\sum_{\xi\in\cX_j}|R_\xi|\sup_{x\in R_\xi}|g(x)|^p\Big)^{1/p}
= \Big\|\sum_{\xi\in\cX_j}M_\xi \ONE_{R_\xi}\Big\|_p
\le c \Big\|\sum_{\eta\in\cX_{j+\ell}}m_\eta^* \ONE_{R_\eta}\Big\|_p\\
&\qquad\quad
\le c \Big\|\cM_s\Big(\sum_{\eta\in\cX_{j+\ell}}m_\eta \ONE_{R_\eta}\Big)\Big\|_p
\le c \Big\|\sum_{\eta\in\cX_{j+\ell}}m_\eta \ONE_{R_\eta}\Big\|_p
\le c\|g\|_p. \qed
\end{align*}


\end{document}